\documentclass[11pt]{article}
\usepackage{latexsym}
                 \usepackage[all]{xy}
                       \usepackage{amscd}
                       \usepackage{amssymb}
                       \usepackage{amsmath}
                       \usepackage{amsthm}
                       \usepackage{enumerate}
                       \usepackage{verbatim}
\def\co{\colon\thinspace}

\def\Der{\mathrm{Der}}
\def\id{\mathrm{id}}
\newtheorem{thm}{Theorem}[section]

\newtheorem{cor}[thm]{Corollary}

\newtheorem{lem}[thm]{Lemma}
\newtheorem{prop}[thm]{Proposition}
\newtheorem{defn}[thm]{Definition}

\newtheorem{prob}[thm]{Problem}

\newtheorem{Example}[thm]{Example}
\newenvironment{ex}{\begin{Example}\rm}{\end{Example}}
\newtheorem{Counterexample}[thm]{Counterexample}

\newtheorem{remark}[thm]{Remark}
\newenvironment{rmk}{\begin{remark}\rm}{\end{remark}}
\newtheorem{Fact}[thm]{Fact}

\newtheorem{Nothing}[thm]{$\!\!\!$}

\newcommand{\be}{\begin{equation}}
\newcommand{\ee}{ \end{equation}}
\newcommand{\ba}{\begin{eqnarray}}
\newcommand{\ea}{\end{eqnarray}}
\newcommand{\ban}{\begin{eqnarray*}}
\newcommand{\ean}{\end{eqnarray*}}

\begin{document}
\abovedisplayskip=6pt plus3pt minus3pt \belowdisplayskip=6pt
plus3pt minus3pt
\title{\bf Obstructions to nonnegative curvature and rational
homotopy theory
\rm}
\author{Igor Belegradek \and Vitali Kapovitch}
\date{}
\maketitle
\begin{abstract}
We establish a link between rational homotopy theory and the
problem which vector bundles admit complete Riemannian metric of
nonnegative sectional curvature. As an application, we show for a
large class of simply-connected nonnegatively curved manifolds
that, if $C$ lies in the class and $T$ is a torus of positive
dimension, then ``most'' vector bundles over $C\times T$ admit no
complete nonnegatively curved metrics.
\end{abstract}

\section{Introduction}
\label{sec: intr-results}
According to the soul theorem of J.~Cheeger and
D.~Gromoll~\cite{CG}, a complete open manifold of nonnegative
sectional curvature is diffeomorphic to the total space of the
normal bundle of a compact totally geodesic submanifold, called
the soul.

A natural problem is to what extent the converse to the soul
theorem holds. In other words, one asks which vector bundles admit
complete nonnegatively curved metrics.  Various aspects of this
problem have been studied
in~\cite{Che, Ri, OW, Yan, GZ, BK1}.
In this paper we only deal with bundles over closed manifolds
diffeomorphic to $C\times T$, where $C$ is simply-connected
and $T$ is a standard torus. By~\cite{CG} any soul has a finite cover
of this form, where $\sec(C)\ge 0$.

Until recently, obstructions
to the existence of metrics with $\sec\ge 0$ on vector bundles
were only known for flat souls~\cite{OW}, which corresponds to
the case when $C$ is a point.
In~\cite{BK1} we produced a variety of
examples of vector bundles
which do not admit complete metrics with
$\sec\ge 0$. For instance, we showed that for
any $C$ and $T$ with $\dim(T)\ge 4$ and $k\ge 2$,  there are
infinitely many rank $k$ vector bundles over $C\times T$
whose total spaces admit no complete metrics with $\sec\ge 0$.
In all these examples $\dim(T)>0$, in fact no obstructions
are known to the existence of complete
metrics with $\sec\ge 0$ on vector bundles
over simply-connected nonnegatively curved manifolds.

We first explain our approach to finding obstructions
in case $\dim(T)>0$.
The main geometric ingredient
is the splitting theorem in~\cite{Wil, BK1},
which says that after passing to a finite cover,
the normal bundle to the soul can be taken,
by a base-preserving diffeomorphism, to the product
$\xi_{C}\times T$ of a vector bundle $\xi_{C}$ over $C$ with
the torus $T$ (see Theorem~\ref{geomsplit}).
Then one is faced with the purely
topological problem of recognizing whether a given vector bundle
over $C\times T$ has this property.
In other words, one needs to study the orbit of
$\xi_{C}\times T$ under the action of the
diffeomorphism group of $C\times T$.
Since vector bundles are rationally classified by
the Euler and Pontrjagin classes, the problem reduces
to analyzing the action of $\mathrm{Diffeo}(C\times T)$
on the rational cohomology algebra $H^*(C\times T,\mathbb Q)$
of $C\times T$.
The ``Taylor expansion'' in $T$-coordinates of
any self-diffeomorphism of $C\times T$ gives rise to a
negative degree derivation of $H^*(C,\mathbb Q)$.
One of the main points of this paper is that
the orbit of $\xi_{C}\times T$ consists of bundles of the
same form, unless there exists a negative degree derivation of
$H^*(C,\mathbb Q)$ that does not vanish on the
Euler or Pontrjagin classes of $\xi_C$.
In particular, if $H^*(C,\mathbb Q)$ has no nonzero negative degree
derivations, the above topological problem gets solved,
which immediately implies that
``most'' bundles over $C\times T$ admit no
complete metric of $\sec\ge 0$.

To state our main results we need the following technical
definition. Given a vector bundle $\xi$ over $C\times T$,
we say that $\xi$ {\it virtually comes from} $C$ if  for some finite
cover $p\co T\to T$, the pullback of $\xi$ by $\id_C\times p$ is
isomorphic to the product $\xi_C \times T$ where $\xi_C$ is a
bundle over $C$.

If $\xi$ virtually comes from $C$, then no known
method can rule out the existence of a complete metric with
$sec\ge 0$ on the total space $E(\xi)$ of $\xi$,
and potentially all such bundles
might be nonnegatively curved.

In this paper
we show that the converse is often true, namely,
under various assumptions on $C$, we show that,
if $\xi$ is a vector bundle over $C\times T$
such that $E(\xi)$ admits a complete metric with $\sec\ge 0$,
then $\xi$ virtually comes from $C$.
This happens for any $C$ if $\xi$ has rank two.

\begin{thm}\label{thm: rk2-split-rigid}
Let $C$ be a closed smooth simply-connected manifold, and $T$
be a torus. Let $\xi$ be a rank two vector bundle over $C\times T$.
If $E(\xi)$ admits a complete metric with $\sec\ge 0$,
then $\xi$ virtually comes from $C$.
\end{thm}

Oriented $\mathbb R^2$-bundles over $C\times T$ are in one-to-one
correspondence via the Euler class with
$H^2(C\times T,\mathbb Z)\cong H^2(C,\mathbb
Z)\oplus H^2(T,\mathbb Z)$. Thus, any oriented
$\mathbb R^2$-bundle $\xi$ over $C\times T$ can be written
uniquely as $c_\xi +t_\xi$ where $c_\xi\in H^2(C,\mathbb Z)$,
$t_\xi\in H^2(T,\mathbb Z)$.
Theorem~\ref{thm: rk2-split-rigid}
implies that if $\sec(E(\xi))\ge 0$, then $t_\xi=0$.

More generally,
it is easy to see that ``most'' vector bundles over $C\times T$ do
not virtually come from $C$, at least when $\dim(T)$ is large
enough (for a precise result,
see~\cite[4.4, 4.6]{BK1} and Lemma~\ref{lem: bundles} below).
In fact, $\xi$ virtually comes from $C$ iff all rational characteristic
classes of $\xi$ lie in the
$H^*(C,\mathbb Q)\otimes H^0(T,\mathbb Q)$-term of the
K\"unneth decomposition $\oplus_{i}H^*(C,\mathbb Q)\otimes H^i(T,\mathbb Q)$
of $H^*(C\times T,\mathbb Q)$.

One of the main sources of examples  of closed manifolds of
nonnegative curvature is given by homogeneous spaces or, more
generally, biquotients of compact Lie groups.
In this case we prove

\begin{thm}\label{biquotients}
Let $C=G//H$ be a simply connected biquotient of
compact  Lie groups such that $H$ is
semi-simple, and let $T$ be a torus.
Let $\xi$ be a vector bundle over $C\times T$ of rank $\le 4$.
If $E(\xi)$ admits a complete metric with $\sec\ge 0$,
then $\xi$ virtually comes from $C$.
\end{thm}

Let $\mathcal H$ be the class of simply-connected finite CW-complexes
whose rational cohomology algebras
have no nonzero derivations of negative degree.

\begin{thm}\label{H-class}
Let $C\in\mathcal H$ be a closed smooth manifold,
and $T$ be a torus.
If $\xi$ is a vector bundle over $C\times T$
such that $E(\xi)$ admits a complete metric with $\sec\ge 0$,
then $\xi$ virtually comes from $C$.
\end{thm}

For example, $\mathcal H$ contains any compact simply-connected
K\"ahler manifold~\cite{Mei2}, and any compact homogeneous space
$G/H$ such that $G$ is a compact connected Lie group and $H$ is a
closed subgroup with $\mathrm{rank}(H)=\mathrm{rank}(G)$~\cite{ST}.
It was proved
in~\cite{Mar} that the total space of a fibration belongs to
$\mathcal H$ provided so do  the base and the fiber.

Recall that a finite
simply-connected cell complex $C$ is called {\it elliptic} if all
but finitely many homotopy groups of $C$ are finite. If $C$ is
elliptic, then $C$ has
nonnegative Euler characteristic, and the sum of the Betti
numbers of $C$ is $\le 2^m$ where $m$ is the cohomological dimension of
$C$~\cite{Fel}.

Any compact simply connected homogeneous space or biquotient is
elliptic, and more
generally, all {\it known} closed simply connected nonnegatively
curved manifolds are elliptic. In fact, it is conjectured~\cite{GH}
that any closed simply connected nonnegatively curved manifold is elliptic. If
true, the conjecture would imply a classical conjecture of
Chern-Hopf that nonnegatively curved manifolds have nonnegative
Euler characteristic, and a conjecture of Gromov that the sum of the
Betti numbers of a compact nonnegatively curved $m$-manifold is
$\le 2^m$.

Halperin conjectured that any elliptic space $C$ of positive Euler
characteristic belongs to $\mathcal H$. This conjecture, which
is considered one of the central problems in rational homotopy theory,
has been confirmed in several important cases~\cite[page 516]{FHT}.
Note that if the above conjectures are true, then
$\mathcal H$ contains any simply-connected compact nonnegatively
curved manifold of positive Euler characteristic.

We refer to the body of the paper for other results similar to
Theorems~\ref{thm: rk2-split-rigid}--\ref{H-class}.
In particular, in
Section~\ref{sec: spherebundles} we establish analogs of
Theorem~\ref{H-class} for $C$'s that belong to
several classes of sphere bundles.
Note that sphere bundles over closed
nonnegatively curved manifolds are potentially a good source of compact
manifolds with $\sec\ge 0$, because
the unit sphere bundle of the normal bundle to the soul is nonnegatively
curved~\cite{GW}.
Also in Section~\ref{sec: pos curv}, we prove an analog of
Theorem~\ref{H-class} where $C$ is any {\it currently known}
simply connected positively curved manifold.

Not every nonnegatively curved vector bundle over $C\times T$
virtually comes from $C$, even though finding an explicit
counterexample is surprisingly difficult. In fact, the conclusion
of Theorem~\ref{biquotients} fails already for rank six bundles
over homogeneous spaces.

\begin{thm}\label{mainex}\it\!
Let $C=SU(6)/(SU(3)\times SU(3))$ and $\dim(T)\ge 2$,
then there exists a rank six vector
bundle $\xi$ over $C\times T$ which does not virtually
come from $C$, but $E(\xi)$ admits a complete metric of $\sec\ge 0$
such that the zero section is a soul.\rm
\end{thm}

To prove the above theorem, we find a nonnegatively curved vector
bundle $\xi_C$ over $C$ with the zero section being a soul, and a
negative degree derivation $D$ of $H^*(C,\mathbb Q)$ that is
induced by a derivation of the minimal model of $C$, and
furthermore such that $D$ does not vanish on the Euler class of
$\xi_C$, but vanishes on the Pontrjagin classes of the tangent
bundle of $C$. Finding such $\xi_C$ and $D$ is not easy, and what
makes it work here is some very special properties of the minimal
model of $C$. Incidentally, $SU(6)/(SU(3)\times SU(3))$ is one of
the simplest non-formal homogeneous spaces.

Now since $D$ is induced by a derivation of the minimal model, a
multiple of $D$ can be ``integrated'' to a self-homotopy
equivalence $f$ of $C\times T$. Furthermore, $f$ preserves the
Pontrjagin classes of the tangent bundle of $C\times T$, because
$D$ vanishes on the Pontrjagin classes of $TC$. Then by a
surgery-theoretic argument, some iterated power of $f$ is
homotopic to a diffeomorphism. Finally, since $D$ does not vanish
on the Euler class of $\xi_C$, the $f$-pullback of $\xi_C\times T$
does not virtually come from $C$, yet it carries the pullback
metric of $\sec\ge 0$ with zero section being a soul.

\

\bf{Structure of the paper.\ }\rm

Section~\ref{sec: notations} is a list of notations and
conventions we use throughout the paper.

In Section~\ref{sec: split criterion} we introduce
a purely topological property of a triple $(C,T,k)$,
which we call splitting rigidity. As a link to nonnegative
curvature, we show that if $(C,T,k)$ is splitting rigid,
and $\xi$ is a rank $k$ vector bundle over $C\times T$
such that $E(\xi)$ admits a complete metric with $\sec\ge 0$,
then $\xi$ virtually comes from $C$.

In Section~\ref{sec: taylor-in-cohom} we relate splitting rigidity
to the absence of negative degree derivations of the cohomology
algebra of $C$.
Sections~\ref{sec: rk2} --~\ref{sec: pos curv} are devoted
to applications, in particular, here we prove the results stated in
Section~\ref{sec: intr-results}, as well as splitting rigidity for
certain sphere bundles and for all known positively curved manifolds.

Section~\ref{sec: der-in-min} is an in depth study of splitting
rigidity. In particular,  we prove that if $k$ is sufficiently large,
then
splitting rigidity can be
expressed in rational homotopy-theoretic terms.
Section~\ref{sec: non-splil-rigid-ex} contains the proof of
Theorem~\ref{mainex}.

In Section~\ref{sec: ex} we show by example that if $k$ is small,
then changing $C$ within its homotopy type
may turn a splitting rigid triple into a non-splitting rigid one.
In Section~\ref{sec: souls-sections} we obtains stronger obstructions
to nonnegative curvature on a vector bundle under the assumption that
the zero section is a soul.

In Section~\ref{sec: prob} we pose and discuss several open problems.
The appendix contains a surgery-theoretic lemma and an existence
result for vector bundles with prescribed Euler and Pontrjagin classes.

Much of the paper can be read without any knowledge of rational
homotopy theory. In fact, rational homotopy background is only needed
for Sections~\ref{sec: biquot},~\ref{sec: pos curv},~\ref
{sec: der-in-min},~\ref {sec: non-splil-rigid-ex} and~\ref{sec: prob}.

\

\bf{Acknowledgments.\ }\rm

The first author is grateful to McMaster University
and California Institute of Technology
for support and excellent working conditions.

It is our pleasure to thank Gregory Lupton and Samuel Smith for
insightful discussions on rational homotopy theory, Stefan
Papadima for Lemma~\ref{vanish}, Ian Hambleton for Lemma~\ref{lem:
surgery}, Toshihiro Yamaguchi for Example~\ref{ex:Yam}, Alexander
Givental for incisive comments on deformation theory, and Burkhard
Wilking and Wolfgang Ziller for countless discussions and insights
related to this work. We are grateful to the referee for helpful
advice on the exposition.
As always, the authors are solely
responsible for possible mistakes.

The present paper grew out of our earlier preprint~\cite{BK2} written in
the summer of 2000. In~\cite{BK2} we proved much weaker results,
say, Theorem~\ref{H-class} is stated there as an open question.
Most of the results of the present paper were obtained
in May and early June of 2001, and reported by the first author during the
Oberwolfach geometry meeting on June 12, 2001. On June 29, 2001
we received a preprint by Jianzhong Pan where he independently proves
Theorem~\ref{H-class} in response to our question in~\cite{BK2}.
We are very grateful to Berhard Hanke for noticing an error in the original proof of Lemma~\ref {lem: bundles} and to Berhard Hanke and Neil Strickland for their help in fixing the error.
\section{Notations and conventions}
\label{sec: notations}
Unless stated otherwise, all (co)homology groups have rational
coefficients, all characteristic classes are over rationals,
all manifolds and vector bundles are smooth.

Given a  cell complex $X$,
define $Char(X,k)$ to be the subspace of $H^*(X)$ equal to
$\oplus_{i=1}^{m} H^{4i}(X)$ if $k=2m+1$, and
equal to $(\oplus_{i=1}^{m-1}H^{4i}(X))\oplus H^{2m}(X)$ if $k=2m$.
If $\xi$ is a (real) oriented rank $k$ vector bundle over $X$,
then in case $k$ is odd, $Char(X,k)$ contains the
Pontrjagin classes $p_1(\xi),\dots, p_m(\xi)$, and
in case $k$ is even,  $Char(X,k)$ contains the
Pontrjagin classes $p_1(\xi),\dots, p_{m-1}(\xi)$, and
the Euler class $e(\xi)$.
The total Pontrjagin class $\sum_{i\ge 0} p_i(\xi)$ is denoted by
$p(\xi)$.

For the product $X\times Y$ of pointed spaces $X$, $Y$,
we denote the projections of $X\times Y$
onto $X$, $Y$ by $\pi_X$, $\pi_Y$. The basepoints
define the inclusions of $X$, $Y$ into $X\times Y$
which we denote by $i_X$, $i_Y$.
For a map $f\co X\times Y\to X^\prime\times Y^\prime$,
we define $f_{XX^\prime}=\pi_{X^\prime}\circ f\circ i_X$,
and $f_{YY^\prime}=\pi_{Y^\prime}\circ f\circ i_Y$.

For the rest of the paper, $C$ stands for a closed, connected,
simply-connected, smooth manifold, and $T$ stands for a torus
of some positive dimension.
We use the K\"unneth isomorphism
$H^*(C\times T)\cong H^*(C)\otimes H^*(T)$ to identify
$\pi_C^*(H^*(C))$ with
the subalgebra $H^*(C)\otimes 1$. We denote the total space
of a vector bundle $\xi$ by $E(\xi)$.

\section{Splitting criterion}
\label{sec: split criterion}
The main geometric ingredient used in this paper is
the following splitting theorem proved in~\cite{BK1}
(cf.~\cite{Wil}).

\begin{thm}\label{geomsplit}
Given a soul $S$ of an open complete nonnegatively curved
manifold $M$, there is a finite cover $p\co\tilde M\to M$,
a soul $\tilde S$ of $\tilde M$ satisfying $p(\tilde S)=S$,
and a diffeomorphism $f\co \tilde S\to C\times T$,
where $C$ is a simply-connected manifold with $\sec(C)\ge 0$
and $T$ is a torus, such that the normal bundle to $\tilde S$
is the $f$-pullback of the bundle $\xi_C\times T$,
where $\xi_C$ is a vector bundle over $C$ whose
total space admits a metric of nonnegative
curvature with the zero section being a soul.\end{thm}

Let $\xi$ be a
vector bundle over $C\times T$. We say that $\xi$ satisfies
($\ast$) if

\hfill $ \begin{array}{l} E(\xi) \text { has a finite cover
diffeomorphic to the
product of T and the total } \\ \text{space of a vector bundle over a closed
simply-connected manifold.}
\end{array}  $ \hfill $(\ast) $

We seek to understand how  assumption ($\ast$) restricts $\xi$.
In particular, we want to find conditions on $C$ ensuring that if
$\xi$ satisfies ($\ast$), then $\xi$ virtually comes from $C$.

\begin{defn}
A triple $(C,T,k)$, where  $k>0$ is an integer, is called {\rm
splitting rigid},
if any rank $k$ vector bundle $\xi$ over $C\times T$ that satisfies
($\ast$) virtually comes from $C$.
\end{defn}

By Theorem~\ref{geomsplit}, if $\sec(E(\xi))\ge 0$,
then $\xi$ satisfies ($\ast$), so we get:

\begin{prop} \label{split-rigid-implies-virt-split}
If $(C,T,k)$ is splitting rigid, and $\xi$ is a rank $k$
vector bundle over $C\times T$ such that $E(\xi)$
has a complete metric with $\sec\ge 0$, then $\xi$ virtually comes
from $C$.
\end{prop}

Thus if $(C,T,k)$ is splitting rigid, then the
total spaces of ``most'' rank $k$ vector bundles over
$C\times T$ do not admit complete metrics with $\sec\ge 0$.

As we prove in Section~\ref{sec: der-in-min}, splitting rigidity
can be often expressed in rational homotopy-theoretic terms.
For example, if $k\ge\dim(C)$, then a triple $(C,T,k)$
is splitting rigid if and only if,
for any derivation of the minimal model of $C$ that commutes
with the differential and has degree within $[-\dim(T),0)$,
the induced derivation on $H^*(C)$ vanishes on $Char(C,k)$.

The same statement holds with some other assumptions in place of
$k\ge\dim(C)$, such as ``$p_i(TC)\in Char(C,k)$ for all $i>0$''.
In Section~\ref{sec: ex},
we give an example when the ``only if'' part fails for $k<\dim(C)$.
On the other hand, the ``if'' part is true without any assumptions on $k$.
As a first step towards these results, we prove
the following proposition, whose weak converse is obtained in
Section~\ref{sec: der-in-min}.

\begin{prop} \label{characterization}
If any self-homotopy equivalence
of $C\times T$ maps $Char(C,k)\otimes 1$ to itself,
then $(C,T,k)$ is splitting rigid.
\end{prop}
\begin{proof}
To check that $(C,T,k)$ is splitting rigid, we
need to start with an arbitrary rank $k$ vector bundle $\xi$
over $C\times T$ that satisfies ($\ast$), and prove that
$\xi$ virtually comes from $C$.
Without loss of generality, we can pass to a finite
cover to assume that $E(\xi)$ is the total space of
a vector bundle $\eta$, which is the product of $T$ and a vector
bundle over a closed smooth simply-connected manifold $C^\prime$.
In other words, $\eta$ is the $\pi_{C^\prime}$-pullback of
a vector bundle over $C^\prime$.

Fix base points in $C$, $C^\prime$, $T$ so that
the inclusions $i_C$, $i_C^\prime$, $i_T$ are defined, and
let $B=C\times T$, $S=C^\prime\times T$.
We think of $\xi$ and $\eta$ as two vector bundle structures
on a fixed manifold $N$, and use the zero sections to identify
$B$ and $S$ with smooth submanifolds of $N$, and identify
$\xi$, $\eta$ with the normal bundles to $B$ and $S$.
Note that both $C$ and $C^\prime$ are
homotopy equivalent to the universal cover of $N$.

Let $g\co B\to S$ be the homotopy equivalence induced
by the zero section of $\xi$ followed by the projection of $\eta$.
Note that $\eta$ is orientable, as the pullback of a bundle
over a simply-connected manifold. Fix orientations of $S$ and $\eta$,
which defines an orientation on $E(\eta)=E(\xi)$.
We orient $B$ so that $\deg(g)=1$, which defines an orientation on $\xi$.

To simplify notations, assume that $k=2m$; the case of odd $k$
is similar. Since $\xi$ has rank $k$,
$e(\xi), p_i(\xi)\in Char(B,k)$ for $0<i<m$.
As we remarked in~\cite[4.4]{BK1},
the bundle $\xi$ virtually comes from $C$ iff
$e(\xi)$ and  $p(\xi)$ lie in $H^*(C)\otimes 1\subset H^*(C\times T)$.
Alternatively, since $Char(B,k)\cap (H^*(C)\otimes 1)=Char(C,k)\otimes 1$,
we see that $\xi$ virtually comes from $C$ iff
$e(\xi), p_i(\xi)\in Char(C,k)\otimes 1$ for $0<i<m$.

By Whitehead's theorem, the maps
$g_{CC^\prime}\co C\to C^\prime$, $g_{TT}\co T\to T$
are homotopy equivalences.
Fix their homotopy inverses $g^{-1}_{CC^\prime}$, $g^{-1}_{TT}$.
Consider a self-homotopy equivalence
$h=(g^{-1}_{CC^\prime}\times g^{-1}_{TT})\circ g$ of $C\times T$.
Note that each of the maps $h_{CC}$, $h_{TT}$ is homotopic to
the identity because, say, $h_{CC}$ is equal to
\[\pi_C\circ h\circ i_C\simeq\pi_C\circ
(g^{-1}_{CC^\prime}\times g^{-1}_{TT})\circ g\circ i_C\simeq
g^{-1}_{CC^\prime}\circ \pi_{C^\prime}\circ g\circ i_C\simeq
g^{-1}_{CC^\prime}\circ g_{CC^\prime}\simeq \mathrm{id}_C\]
Now it is routine to check that $g^*$
maps $Char(C^\prime,k)\otimes 1$ into $Char(C,k)\otimes 1$
if and only if $h^*$ maps $Char(C,k)\otimes 1$ to itself.

By assumption $h^*$ maps $Char(C,k)\otimes 1$ to itself,
so $g^*(Char(C^\prime,k)\otimes 1)=Char(C,k)\otimes 1$.
It was observed in~\cite[section~3]{BK1} that
$g^*$ maps $e(\eta)$ to $e(\xi)$. Thus $e(\xi)\in Char(C,k)\otimes
1$, as needed, and
it remains to show that $Char(C^{\prime},k)\otimes 1$ also contains
$p_i(\xi)$ for $0<i<m$.

Since $g$, viewed as a map $B\to N$,
is homotopic to the inclusion $B\hookrightarrow N$, we
have that $TN|_{B}\cong g^{\#}(TN|_S)$.
By the Whitney sum formula
\[
p(TB)p(\xi)=p(TN|_B)=p(g^{\#}(TN|_{S}))=p(g^{\#}(TS\oplus
\eta))=g^{*}(p(TS))g^{*}(p(\eta))
\]
Since $T$ is parallelizable, $p(TB)=p(TC)\otimes 1\in H^*(C)\otimes 1$
and $p(TS)=p(TC^\prime)\otimes 1\in H^*(C^{\prime})\otimes 1$.
Since $p(TB)$ is a unit in $H^*(C)\otimes 1$, we can write
$p(TB)^{-1}=\sum_j a_j\otimes 1$ for some $a_j\in H^*(C)$, and
\begin{equation*}
p(\xi)=p(TB)^{-1}g^{*}(p(TS))g^{*}(p(\eta))=\left(\sum_{j\ge 0}
a_j\otimes 1\right)
\left(\sum_{l\ge 0}g^{*}(p_l(TC)\otimes 1)\right)
\left(\sum_{n\ge 0}g^{*}(p_n(\eta))\right)
\end{equation*}
Also $p_m(\eta)\in H^{*}(C^{\prime})\otimes 1$, because $\eta$ is
a product of torus and a bundle over $C^{\prime}$. By definition
of ``Char'', this means that $p_l(TC^\prime)\otimes 1$ and
$p_n(\eta)$ lie in $Char(C^{\prime},k)\otimes 1$ for any
$0<l,n<m$, and therefore, $g^{*}(p_l(TC^\prime)\otimes 1))$,
$g^{*}(p_n(\eta))\in Char(C,k)\otimes 1$ for any $0<l,n< m$. The
above formula now implies that $p_i(\xi)\in Char(C,k)\otimes 1$
for $0<i<m$, as promised. This completes the proof that $\xi$
virtually comes from $C$.
\end{proof}
\begin{rmk}\label{hom-to-id} It is easy to see that the assumption of
Proposition~\ref{characterization} that $Char(C,k)\otimes 1$ is
invariant under any self-homotopy equivalence of $C\times T$ is
equivalent to the
formally weaker assumption that $Char(C,k)\otimes 1$ is
invariant under any self-homotopy equivalence of $C\times T$ satisfying
$h_{CC}\sim \mathrm{id}_{C}$, $h_{TT}\sim \mathrm{id}_{T}$.
\end{rmk}
\begin{rmk}
Proposition~\ref{characterization} implies that
if $Char(C,k)\otimes 1$ is invariant under
any graded algebra automorphism of
$H^*(C\times T)$, then $(C,T,k)$ is splitting rigid, and
this is how we establish splitting rigidity in this paper,
except for one example in Section~\ref{sec: ex}
where deeper manifold topology gets involved.
\end{rmk}
\begin{ex} \label{easy ex of split rigidity}
Of course,~\ref{characterization} applies if $Char(C,k)=0$.
Thus, $(S^{2m+1}, T, k)$ is splitting rigid for any $k$, $T$.
\end{ex}
In a special case $e(\xi)\in H^*(C)\otimes 1$, the very same proof of
Proposition~\ref{characterization} yields
the following stronger statement, in which $\epsilon^{m}$
denotes the  trivial rank $m$ bundle over $C\times T$.
\begin{prop}\label{prop: stabilizing}
Let $\xi$ be a rank $k$ vector bundle over $C\times T$
with $e(\xi)\in H^*(C)\otimes 1$ such that $\xi\oplus\epsilon^{m}$
satisfies ($\ast$) for some $m\ge 0$.
If any self-homotopy equivalence
of $C\times T$ maps $Char(C,k)\otimes 1$ to itself, then
$\xi$ virtually comes from $C$.
\end{prop}

The assumption on $e(\xi)$ cannot be dropped: say
if $\xi$ is a rank two bundle over the $2$-torus with
$e(\xi)\neq 0$, then $\xi\oplus\epsilon^1$ becomes trivial
in a finite cover, but $\xi$ does not. Of course, if $e(\xi)=0$
(which for example is always true if $k$ is odd), then $e(\xi)\in
H^*(C)\otimes 1$.

\section{Taylor expansion in cohomology}
\label{sec: taylor-in-cohom}

Let $A=\oplus_{p}A_p$ be a graded $\mathbb Q$-algebra.
If $a\in A_p$, we refer to $p$
as the {\it degree} of $a$ and denote it by $|a|$.
In this paper we only consider graded commutative algebras
with an identity element $1\in A_0$, and such that
$A_p=0$ for $p<0$.

Let $B$ be a subalgebra of $A$, and let $n\in\mathbb Z$.
A {\it degree $n$ derivation of $B$ with values in $A$}
is a linear map $D\co B\to A$ such that
if $a\in A_p$, then $|D(a)|=p+n$, and
$D(ab)=D(a)b+aD(b)(-1)^{np}$ for any $a\in A_p$, $b\in A$.
We refer to $n$ as the degree of $D$, and denote it by $|D|$.
If $B=A$, we just say that $D$ is a {\it degree $n$
derivation} of $A$.

Let $\mathrm{Der}_n(B,A)$
be the $\mathbb Q$-vector space of
degree $n$ derivation of $B$ with values in $A$, and
write $\mathrm{Der}_n(A)$ for $\mathrm{Der}_n(A,A)$.
Let $\mathrm{Der}_-(A)=\oplus_{n<0}\mathrm{Der}_n(A)$.
We refer to derivations of $A$ of negative degree
(i.e. to elements of $\mathrm{Der}_-(A)$) as
{\it negative} derivations of $A$.

The cohomology algebra $H^*(T)$ of $T$ is an exterior algebra
on degree one generators $x_j$ with $j=1,\dots ,\dim(T)$.
The $\mathbb Q$-vector space $H^*(T)$ has an
obvious basis $\{t_i\}$, $i=0,\dots, n$ of square-free monomials
in variables $x_j$ where $n=2^{\dim(T)}-1$.
We order $\{t_i\}$ lexicographically so that
$t_0=1$, $t_j=x_j$ for $j=1,\dots ,\dim(T)$,
and $t_n=x_1\ldots x_{\dim(T)}$. Thus, $t_i^2=0$ for $i>0$.
We write $H^*(T)=\oplus_{i}\mathbb Qt_i$.
Then $H^*(C\times T)=H^*(C)\otimes H^*(T)=
\oplus_{i} H^*(C)\otimes\mathbb Qt_i$
is a free $H^*(C)$-module.

Let $h$ be a self-homotopy equivalence of $C\times T$. Since
$H^*(C\times T)$ is free $H^*(C)$-module with basis $\{t_i\}$,
given $a\in H^*(C)$, there is a unique sequence of elements
$\frac{\partial h^*}{\partial t_i}(a)\in H^*(C)$ such that
$h^*(a\otimes 1)=\sum_{i}(1\otimes t_i) (\frac{\partial
h^*}{\partial t_i}(a)\otimes 1)$. We think of $\frac{\partial
h^*}{\partial t_i}$ as $\mathbb Q$-linear self-maps of $H^*(C)$.
Informally, it is useful to interpret the above formula as a
Taylor expansion of $h^*$ at $a\otimes 1$.

Since $h^*$ is an algebra isomorphism, the maps $\frac{\partial
h^*}{\partial t_i}$ satisfy certain recursive identities, obtained
from $h^*(ab)=h^*(a)h^*(b)$ by collecting the terms next to
$1\otimes t_i$'s. For example, $\frac{\partial h^*}{\partial t_0}$
is an algebra isomorphism of $H^*(C)$, and
\begin{equation}\label{e:hom}
\frac{\partial h^*}{\partial t_1}(ab)= \frac{\partial
h^*}{\partial t_1}(a) \frac{\partial h^*}{\partial
t_0}(b)+(-1)^{-|a|} \frac{\partial h^*}{\partial t_0}(a)
\frac{\partial h^*}{\partial t_1}(b).
\end{equation}
If $h_{CC}\sim\mathrm{id}_C$, which can always be arranged
in our case by Remark~\ref{hom-to-id},
then $\frac{\partial h^*}{\partial t_0}(a)=a$, and
therefore, $\frac{\partial h^*}{\partial t_1}$
is a derivation of $H^*(C)$ of degree $-1$.

Let $d=\dim T$. Suppose all degree $-1$ partial derivatives
$\frac{\partial h^*}{\partial t_1}, \ldots, \frac{\partial
h^*}{\partial t_{d} }$ vanish. Then we claim that $\frac{\partial
h^*}{\partial t_{d+1} }$ is a derivation.
Indeed, since $h^*$ is a homomorphism we get

\begin{gather}
ab\otimes 1+\frac{\partial h^*}{\partial t_{d+1} }(ab)\otimes
t_{d+1}+\text{higher order terms}=h^*(ab\otimes 1)=h^*(a\otimes
1)h^*(b\otimes 1)=\nonumber\\
   =(a\otimes 1+\frac{\partial h^*}{\partial
t_{d+1} }(a)\otimes t_{d+1}+\text{h. o. terms})(b\otimes
1+\frac{\partial h^*}{\partial t_{d+1} }(b)\otimes
t_{d+1}+\text{h. o. terms})\label{e:der}\\
   =ab\otimes 1+[\frac{\partial
h^*}{\partial t_{d+1} }(a)b +a\frac{\partial h^*}{\partial t_{d+1}
}(b)]\otimes t_{d+1}+\text{h. o. terms}\nonumber
\end{gather}
which proves our assertion.
Similarly, if $\frac{\partial h^*}{\partial t_i}=0$ for $0<i<k$,
then $\frac{\partial h^*}{\partial t_k}$ is a derivation
of degree $-|t_k|$.

More generally, if $h_{CC}$ is not homotopic to $\mathrm{id}_C$
then $\frac{\partial h^*}{\partial t_1}\circ
(\frac{\partial h^*}{\partial t_0})^{-1}$ is a
derivation of $H^*(C)$ of degree $-1$,
and if $\frac{\partial h^*}{\partial t_i}=0$ for $0<i<k$,
then $\frac{\partial h^*}{\partial t_k}\circ
(\frac{\partial h^*}{\partial t_0})^{-1}$ is a
derivation of $H^*(C)$ of degree $-|t_k|$.

Thus, if $H^*(C)$ has no nonzero negative derivations, then
$h^*(a\otimes 1)=a\otimes 1$ for any $a\in H^*(C)$,
and any self-homotopy equivalence $h$ of $C\times T$. Thus,
$(C,T,k)$ is splitting rigid for any $T$, $k$. Combining with
Propositions~\ref{split-rigid-implies-virt-split},~\ref{characterization},
we deduce Theorem~\ref{H-class}.
We actually need the following stronger statement.

\begin{prop}\label{prop: no-der-split-rigid}
If every negative  derivation of $H^*(C)$
vanishes on $Char(C,k)$, then $(C,T,k)$ is splitting rigid
for any $T$.
\end{prop}
\begin{proof}
Let $h$ be a self-homotopy equivalence of $C\times T$ with
$h_{CC}\sim\mathrm{id}_C$, so that
$\frac{\partial h^*}{\partial t_0}(a)=a$, and
$\frac{\partial h^*}{\partial t_1}$
is a derivation of $H^*(C)$ of degree $-1$.
By Proposition~\ref{characterization}, Remark~\ref{hom-to-id},
it suffices to show that $h^*(b\otimes 1)=b\otimes 1$,
for all $b\in Char(C,k)$.

Let $\phi_1$ be a self-map of $H^*(C\times T)$ defined by
$\phi_1(a\otimes t)= a\otimes t-(1\otimes t_1) (\frac{\partial
h^*}{\partial t_1}(a)\otimes t)$ and for $t\in H^*(T)$, $a\in
H^*(C)$. The fact that $\frac{\partial h^*}{\partial t_1}$ is a
derivation and $t_1^2=0$  implies that $\phi_1$ is a homomorphism
(cf. (\ref{e:der}) above). It is also easy to check that the map
$a\otimes t\mapsto a\otimes t+(1\otimes t_1) (\frac{\partial
h^*}{\partial t_1}(a)\otimes t)$ is the inverse to $\phi_1$ and
therefore $\phi_1$ is an automorphism of $H^*(C\times T)$.

   Then
\[\phi_1\circ h^*(a\otimes 1)= a\otimes 1+ \sum_{i\ge 2}(1\otimes
t_i) (\phi_1^i(a)\otimes 1),\] where $\phi_1^i$ are linear
self-maps of $H^*(C)$.

Now $\phi_1^2$ is a derivation of $H^*(C)$,
so the formulas
$\phi_2(a\otimes 1)= a\otimes 1-(1\otimes t_2)
(\phi_1^2(a)\otimes 1)$
and $\phi_2(1\otimes t)=1\otimes t$ for $t\in H^*(T)$,
$a\in H^*(C)$ define an automorphism $\phi_2$ of $H^*(C\times T)$.
Then \[\phi_2\circ \phi_1\circ h^*(a\otimes 1)=
(a\otimes 1)+ \sum_{i\ge 3}(1\otimes t_i)
(\phi_2^i(a)\otimes 1),\]
where $\phi_2^i$ are
linear self-maps of $H^*(C)$, and $\phi_2^3$ is a derivation.
Continuing in this fashion, we get automorphisms
$\phi_k$ with $\phi_k(a\otimes 1)= a\otimes 1-(1\otimes t_i)
(\phi_{k-1}^k(a)\otimes 1)$ where $\phi_{k-1}^k$ is a derivation
of $H^*(C)$, and such that
$\phi_{n}\circ\dots\circ \phi_1\circ h^*(a\otimes 1)=a\otimes 1$.

Thus $h^*(a\otimes 1)=
\phi_{1}^{-1}\circ\dots\circ \phi_n^{-1}(a\otimes 1)$. Also
$\phi_k^{-1}(a\otimes 1)=a\otimes 1+(1\otimes t_i)
(\phi_{k-1}^k(a)\otimes 1)$.
Now if $b\in Char(C,k)$, then by assumption $\phi_{k-1}^k(b)=0$
so that $\phi_k^{-1}(b\otimes 1)=b\otimes 1$. Thus,
$h^*(b\otimes 1)=b\otimes 1$ as desired.
\end{proof}

\section{Splitting rigidity for rank two bundles}
\label{sec: rk2}
The following well-known lemma is the key ingredient in
the proof of Theorem~\ref{thm: rk2-split-rigid}.
\begin{lem}\label{vanish}
Let $A\subset B$ be a finite dimensional  subalgebra of a
commutative graded $\mathbb Q$-algebra $B$
satisfying  $B_0\cong\mathbb Q$. Let $D$ be a derivation of $A$ of
degree $-2n<0$. Then $D$ vanishes on $A_{2n}$.
\end{lem}
\begin{proof}
Let $a\in A^{2n}$. Since $B_0\cong\mathbb Q$ and $|D|=-2n$,
$D(a)\in B_0$ is a rational multiple of $1$.
Choose a positive integer $m$ such that $a^m=0$ but
$a^{m-1}\ne 0$ (which exists since $A$ is finite-dimensional).
Since $|a|$ is even, we get $0=D(a^m)=ma^{m-1}D(a)$ so that $D(a)=0$.
\end{proof}
\begin{proof}[Proof of Theorem~\ref{thm: rk2-split-rigid}]
By Proposition~\ref{split-rigid-implies-virt-split}, it is enough to
show that $(C,T,2)$ is splitting rigid.
Then by Proposition~\ref{prop: no-der-split-rigid}, it suffices to show
that any negative  derivation of $H^*(C)$ vanishes on $H^2(C)$.
Since $C$ is simply connected, negative  derivations of degree
$\ne -2$ automatically vanish on $H^2(C)$ for degree reasons.
Since $H^*(C)$ is finite dimensional,
Lemma~\ref{vanish} implies that all derivations of degree $-2$
vanish on $H^2(C)$ as well.
\end{proof}

\section{Splitting rigidity for sphere bundles}
\label{sec: spherebundles}
In this section we establish splitting rigidity for various classes
of sphere bundles. We need the following standard lemma.
\begin{lem}\label{decomp}
Let $A, B, C$ be graded commutative $\mathbb{Q}$-algebras such that
$A$ is finite dimensional, $B$ is a subalgebra of $C$, and
the algebras $A\otimes B$ and $C$ are isomorphic as $B$-modules.
Let $Der(C)|_B$ be the image of the restriction map
$Der(C)\to Der(B,C)$. Then $Der(C)|_B$ and
$A\otimes Der(B)$ are isomorphic as $\mathbb Q$-vector spaces.
\end{lem}
\begin{proof}[Sketch of the proof]
Let $\{a_i\}$ be a basis of the vector space $A$.
Given a derivation $D\in \mathrm{Der}(C)|_{B}$,
define  linear self-maps $D_i^B$ of $B$ by
$D(1\otimes b)=\sum_i a_i\otimes D_i^B(b)$.
It is routine to check that
$D_i^B\in\mathrm{Der}(B)$, and the correspondence
$D\to\sum_i a_i\otimes D_i^B$
gives the promised isomorphism.
\end{proof}
\begin{rmk} If in the above proof $D\in \mathrm{Der}_{-}(C)$,
then $D_i^B\in\mathrm{Der}_{-}(B)$ for all $i$.
\end{rmk}
\begin{thm}\label{generalproducts}
Let $T$ be a torus, and $C$ be a closed, simply-connected, smooth
manifold which is the total space of a sphere bundle
with zero Euler class and base $B\in \mathcal{H}$
satisfying $H^{odd}(B,\mathbb Q)=0$.
Then if $\xi$ is a
vector bundle over $C\times T$ such that $E(\xi)$ admits a
complete metric with $\sec\ge 0$, then $\xi$ virtually comes from
$C$.
\end{thm}
\begin{proof}
If $l$ is even, both the fiber and the
base of the $S^l$-fibration $p\co C\to B$ belong to $\mathcal{H}$,
and hence $C\in\mathcal{H}$ by~\cite{Mar}.
Thus, we are done by Theorem~\ref{H-class}.

Suppose that $l$ is odd. Since the Euler class of $C\to B$ is
trivial, the Serre spectral sequence of the fibration collapses at
the $E_2$-term, so $H^*(C)$ and $H^*(B)\otimes H^*(S^{l})$ are
isomorphic, as $H^*(B)$-modules.
Since $H^{odd}(B)=0$, this implies that $p\co H^{even}(B)\to
H^{even}(C)$ is an isomorphism. Note that under the above
identification of $H^*(C)$ and $H^*(B)\otimes H^*(S^{l})$,
$p^*(H^*(B))=H^{ev}(C)$ corresponds to $H^*(B)\otimes 1$.

By Proposition~\ref{prop: no-der-split-rigid}, to
prove splitting rigidity  of $(C,T,k)$ for any $T,k$, it is
enough to show that all negative  derivations of $H^*(C)$
vanish on $H^{ev}(C)$.
Let $D\in \mathrm{Der}_{-}(H^*(C))$ be a negative derivation.
By Lemma~\ref{decomp} we can write
$D|_{H^{ev}(C)}$ as $aD_1+D_{2}$ where $a\in H^l(S^{l})$ and
$D_1,D_{2}\in\mathrm{Der}_{-}(H^*(B))$.
Now $B\in\mathcal H$ implies that $D_1=D_{2}=0$, and
hence $D|_{H^{ev}(C)}=0$, as needed.
\end{proof}
In case $\dim(S^{l})\ge\dim(B)$, the Euler class is automatically zero,
therefore we have
\begin{cor} Let $S^{l} \to C\to B$ be a sphere bundle over
$B\in \mathcal{H}$ with $H^{odd}(B)=0$ where $l\ge dim(B)$.
Then $(C,T,k)$ is splitting rigid for any $T,k$.
\end{cor}
\begin{rmk}
If $B$ is elliptic, then $H^{odd}(B,\mathbb Q)=0$
is equivalent  to $\chi(B)>0$~\cite[Proposition 32.10]{FHT}, where
$\chi $ is the Euler characteristic.
In particular, all homogeneous spaces $G/H$ with
$\mathrm{rank}(H)=\mathrm{rank}(G)$ have $H^{odd}(G/H,\mathbb Q)=0$.
As we mentioned in the introduction, conjecturally,
any nonnegatively curved manifold $B$ of positive Euler
characteristic has $H^{odd}(B,\mathbb Q)=0$ and lies in $\mathcal H$.
\end{rmk}

Let $\mathcal H(n)$ be the class of simply-connected
finite CW-complexes whose rational cohomology algebras are
generated in dimension $n$.
By Lemma~\ref{vanish}, $\mathcal H(2n)\subset \mathcal H$ for any $n$.
Notice that this implies that any $B\in \mathcal H(2n)$ satisfies the
assumptions of Theorem~\ref{generalproducts}.
This also implies that the class $\mathcal H(2n)$ is closed
under fibrations because if $F\to E\to B$ is a fibration with the
fiber in $\mathcal H$, then by~\cite{Mar} there is an
$H^*(B)$-module isomorphism $H^*(E)\cong H^*(F)\otimes H^*(B)$.
A simple Mayer-Vietoris arguments implies that if closed simply-connected
manifolds $C$, $C^\prime$ belong to $\mathcal H(2n)$, then so does
the connected sum $C\# C^\prime$.
\begin{ex}
Besides being closed under fibrations and connected sums, the
class $\mathcal H(2n)$ contains the following nonnegatively curved
manifolds: \newline
$\bullet$ $\mathcal H(2)$ contains
$S^2$, $CP^n$, $CP^n\# CP^n$,
$CP^n\#\overline{CP^n}$~\cite{Che}, all nontrivial $S^2$-bundles over
$S^4$~\cite{GZ}, 
biquotients $G//T$ of a compact connected Lie
group $G$ by a maximal torus $T$~\cite{Sin},
$SO(2n+1)/(SO(2n-1)\times SO(2))$ which is a homology $CP^{2n-1}$~\cite{MZ},
projectivized tangent bundle to $CP^n$~\cite{Wil2},
and the
exceptional space $G_2/U(2)$ which is a homology $CP^5$~\cite{MZ};\newline
$\bullet$ $\mathcal H(4)$ contains $S^{4}$, $HP^n$,
$HP^n\# HP^n$,
$HP^n\# \overline{HP^n}$~\cite{Che},
$G_2/SO(4)$ which is a homology $HP^2$, projectivized tangent bundle to
$HP^n$~\cite{Wil2}, and $Sp(n)/K$ where $K$ is the product of $n$
copies of $Sp(1)$; \newline
$\bullet$ $\mathcal H(8)$ contains $S^8$,
$F_4/Spin(9)$, and $F_4/Spin(8)$.
\end{ex}

\begin{thm}\label{spherebundles}
If $S^{l}\to C\to B$ is a sphere bundle such that
$B\in\mathcal H(2n)$ for some $n$,
then $(C,T,k)$ is splitting rigid for any $T,k$.
\end{thm}
\begin{proof}
If $l$ is even, then $S^l, B\in\mathcal H$
so that $C\in\mathcal{H}$ by~\cite{Mar}, and
we are done by Theorem~\ref{H-class}.
If the Euler class of $C\to B$ vanishes, then the
result follows by Theorem~\ref{generalproducts}.
Thus, we can assume that $l$ is odd and
$p\co C\to B$ has nonzero Euler class.

It follows from the Gysin sequence
that $p^*\co H^{ev}(B)\to H^{ev}(C)$ is onto, so
$H^{ev}(C)$ is generated in dimension $2n$.
Again, we see from the Gysin sequence that $H^i(C)=0$ for $0<i<2n$.
Hence, Lemma~\ref{vanish}
implies that $\mathrm{Der}_{-}(H^{ev}(C),H^*(C))=0$.
So by  proposition~\ref{prop: no-der-split-rigid}, $(C,T,k)$ is
splitting rigid for any $T,k$.
\end{proof}

\section{Splitting rigidity for biquotients}
\label{sec: biquot}
Let $G$ be a compact Lie group and $H\leq G\times G$ be a compact
subgroup. Then $H$ acts on $G$ on the left by the formula
$(h_1,h_2)g=h_1gh_2^{-1}$. The orbit space of this action is
called a {\it biquotient} of $G$ by $H$ and denoted by $G//H$. If
the action of $H$ on $G$ is free, then $G//H$ is a manifold. This
is the only case we consider in this paper. In the special case
when $H$ has the form $K_{1}\times K_{2}$ where $K_{1}\subset
G\times 1\subset G\times G$ and $K_{2}\subset 1\times G\subset
G\times G$ we will sometimes write $K_{1}\backslash G/K_{2}$
instead of $G//(K_{1}\times K_{2})$. For any biinvariant metric on
$G$ the above action of $H$ is isometric, and therefore, $G//H$
can be equipped with a submersion metric, which by O'Neill's
submersion formula is nonnegatively curved. Thus, biquotients form
a large class of examples of nonnegatively curved manifolds.

As was observed by Eschenburg~\cite{Esc2}, any
biquotient $G//H$
is diffeomorphic to a biquotient of $G\times G$ by $G\times H$
written as $\Delta G\backslash G\times G/H$, where
$\Delta G$ stands
for the diagonal embedding of $G$ into $G\times G$.
Let $p\co G//H \to B_H$ be the classifying map of the principle
$H$ bundle $H\to G\to G//H$. Then it is easy to see (cf.~\cite{Esc2})
that $G\to G//H\to BH$ is a Serre fibration (which need not be
principal!).
Moreover, this fibration fits into the following fibered square
(see~\cite{Esc2} and~\cite{Sin})
\begin{equation}\label{square}
\xymatrix {G//H\ar[r]\ar[d]&B_G\ar[d]\\
B_H\ar[r]&B_{G\times G}}
\end{equation}
where both vertical arrows are fibrations with fiber $G$ and both
horizontal arrows are fibrations with fiber $(G\times G)/H$. In
particular, the fibration $G//H\to B_H$ is the pullback of the
fibration $G\to B_G\to B_{G\times G}$. Following Eschenburg, we
call the fibration $B_G\to B_{G\times G}$ the {\it reference
fibration}.

Next we are going to construct a Sullivan model of the biquotient
$G//H$.

We refer to~\cite[Chapter1]{OT} for a gentle introduction
to rational homotopy theory, and use the textbook~\cite{FHT} as a
comprehensive reference.

Recall that a free DGA $(\Lambda V, d)$ is called {\it pure} if
$V$ is finite-dimensional and $d|_{V^{ev}}=0$ and $d(V^{odd})\subset
V^{ev}$.
It is well-known that homogeneous spaces admit natural pure
Sullivan models given by their Cartan algebras. The next
proposition shows that the same remains true for biquotients.
\begin{prop}\label{pure}
Let $G//H$ be a biquotient. Then it admits a pure Sullivan model.
\end{prop}
\begin{proof}
We begin by constructing the canonical model of the reference
fibration $\phi\co B_G\to B_{G\times G}=B_G\times B_G$. Since this
fibration is induced by the diagonal map $\Delta \co G\to G\times
G$, it follows that $\phi$ is the diagonal embedding
$\Delta_{B_{G}}\co B_G\to
B_G\times B_G$. Consider the map $\phi^*\co H^*(B_G\times B_G)\to
H^*(B_G)$.  It is well-known that $G$ is rationally homotopy
equivalent to $S^{2m_1-1}\times\dots \times S^{2m_n-1}$ and the
minimal model of $B_G$ is isomorphic to $H^*(BG,\mathbb Q)\cong
\mathbb Q[x_1,...,x_n]$ with zero differentials and with $|x_i|=2m_i$.
Similarly the minimal model of $B_G\times B_G$ is isomorphic to
its cohomology ring $B=\mathbb Q[x_1,...,x_n,y_1,...,y_n]$ with
$|x_i|=|y_i|=2m_i$. Thus $\phi^*$ can be viewed as a
DGA-homomorphism of minimal models of $B_G$ and $B_{G\times G}$.

Let us construct a Sullivan model of $\phi^*$.
Since $\phi=\Delta_{B_{G}}$ we compute
that $\phi^*(x_i)=\phi^*(y_i)=x_i$ for all $i=1,\dots, n$.
Consider the relative Sullivan algebra
$(B\otimes\Lambda(q_1,...q_n),d)$ where $dx_i=dy_i=0$ and
$dq_i=x_i-y_i$. Then it is immediate to check that this relative algebra is
a Sullivan model (in fact, a minimal one) of $\phi^*$ with the
quasi-isomorphism $ B\otimes \Lambda(q_1,...q_n)\to H^*(B_G)$
given by $x_i\to x_i$, $y_i\to x_i$, $q_i\to 0$.

By the naturality of models of maps~\cite[page 204, Proposition 15.8]{FHT},
from the fibered square (\ref{square}), we obtain that
a Sullivan model of the map $G//H\to B_H$ can be given by the
pushout of $(B\otimes\Lambda(q_1,...q_n),d)$ via the homomorphism
$f^*\co B\to H^*(B_H)$; i.e. it can be written as
$$(H^*(B_H),0)\otimes_{(B,d)}(B\otimes
\Lambda(q_1,...q_n),d)=(H^*(B_H)\otimes
\Lambda(q_1,...q_n),\bar{d})$$ where $\bar{d}$ is given by
$\bar{d}|_{H^*(B_H)}=0$ and $\bar{d}(q_i)=f^*(x_i-y_i)$. In
particular, $M(G//H)=(H^*(B_H)\otimes \Lambda(q_1,...q_n),\bar{d})$
is a model for $G//H$. Notice that $H^*(B_H)$ is a free
polynomial algebra on a finite number of even-dimensional
generators, and thus the model $M(G//H)$ is pure.
\end{proof}
\begin{rmk} It is easy to see that the minimal model of a pure
Sullivan model is again pure. Therefore, Proposition~\ref{pure}
implies that minimal models of biquotients are pure.
\end{rmk}
\begin{rmk} The pure model $M(G//H)$  constructed in the proof of
Proposition \ref{pure} provides an effective way of computing
rational cohomology of biquotients. We refer to $M(G//H)$
as the {\it Cartan model} of $G//H$.
This method of computing
$H^*(G//H)$ is essentially equivalent to the method developed by
Eschenburg~\cite{Esc} who computed the Serre spectral sequence
of the fibration $G\to G//H\to B_H$. In fact, it is easy to recover
this spectral sequence by introducing the standard bigrading on
the Cartan model $M(G//H)$. Also note that in case when $G//H$ is
an ordinary homogeneous space (i.e. when $H\subset G\times G$ has
the form $H\times 1\subset G\times G$) this model is easily seen to
reduce to the standard Cartan model of $G/H$.
\end{rmk}
\begin{proof}[Proof of Theorem~\ref{biquotients}]
By Theorem~\ref{thm: rk2-split-rigid}, we only have to consider
the case when $k=3$ or $4$. By
Proposition~\ref{split-rigid-implies-virt-split} it is enough to
show that $(C,T,k)$ is splitting rigid for any $T$ and $k= 3$ or
$4$. Since $Char(C,3)=Char(C,4)=\langle H^4(C)\rangle$ and
according to Proposition~\ref{prop: no-der-split-rigid}, to insure
splitting rigidity it is enough to show that all negative
derivations of $H^{*}(G//H)$ vanish on $H^4(G//H)$.

By passing to a finite cover we can assume that both
$G$ and $H$ are connected.
Since $H$ is semisimple and $G//H$ is simply connected,
the long exact sequence of the fibration $H\to G\to G//H$ implies
that $G$ is also semisimple.

Let $(\Lambda(V),d)$ be the minimal model of $G//H$. Since $G//H$ is
$2$-connected, we have that $V^1=V^2=0$.
According to Proposition~\ref{pure}, the model
$(\Lambda(V),d)$ is a pure. By minimality
of $(\Lambda(V),d)$,
we have that $d|_{V_3}=0$ and $V^3\cong H^3(G//H)$.
By the structure theorem for pure DGAs,~\cite[p 141, Prop 3]{Oni2},
this implies that
$\Lambda(V)\cong \Lambda V^3\otimes \Lambda(\hat{V})$ and
$H^*(\Lambda(V))\cong  \Lambda V^3 \otimes H^*(\Lambda(\hat{V}))$
for some differential subalgebra $\Lambda(\hat{V})\subseteq\Lambda(V)$
such that $\hat{V}^3=0$. Let $A=\Lambda V^3$ and
$B=H^*(\Lambda(\hat{V}))$. By above $B^1=B^2=B^3=0$.
Therefore, by Lemma~\ref{vanish}, all negative derivations of
$B$ vanish on $B^4$. Notice that $H^4(G//H)$
corresponds to $1\otimes B^4$
under the isomorphism $H^*(G//H)\cong A\otimes B$, and hence
applying Lemma~\ref{decomp}, we conclude
that negative derivations of $H^*(G//H)$ vanish on $H^4(G//H)$.
\end{proof}

\begin{rmk}
Using Proposition~\ref{prop: stabilizing}, we get a stronger
version of Theorem~\ref{biquotients}. Namely, if
$\xi$ is a vector bundle over $C\times T$ of rank $3$ or $4$
with $e(\xi)\in H^*(C)\otimes 1$ such that for some $m\ge 0$,
the manifold $E(\xi)\times\mathbb R^m$
admits a complete metric with $\sec\ge 0$, then
$\xi$ virtually comes from $C$.
If $\xi$ has rank $3$, then $e(\xi)=0$, so
the assumption $e(\xi)\in H^*(C)\otimes 1$
is automatically true.
\end{rmk}

\begin{rmk} It would be interesting to see whether
Theorem~\ref{biquotients} remains true if $H$ is not assumed to be
semi-simple. In that case a slight modification of the proof of
Theorem~\ref{biquotients} still shows that derivations of degree
$-4,-3$ and $-1$ vanish on $H^4(G//H)$, and therefore
$(G//H,S^1,k)$ is splitting rigid for $k\le 4$. However, as of
this writing, we are unable to see whether degree $-2$ derivations
have to vanish, and thus the general case remains unclear.
\end{rmk}

\begin{prop}\label{torusbiq}
Let $C=G//H$ be a simply connected biquotient of a compact group $G$
by a torus $H$ satisfying $\mathrm{rank}(H)=\mathrm{rank}(G)-1$. Then
$(C,T,k)$ is splitting rigid for any $T,k$.
\end{prop}
\begin{proof}
First, we show that $H^{even}(G//H)$ is generated in dimension $2$.
Consider the Cartan model
$M(G//H)=(H^*(B_H)\otimes \Lambda(q_1,...q_n),\bar{d})$. It admits a
natural grading  by the wordlength in $q_{i}$'s given by
$M(G//H)_{k}=H^*(B_H)\otimes \Lambda^{k}(q_1,...q_n)$. Since the
differential decreases the wordlength in $q_i$'s by $1$, this grading
induces a natural grading in the cohomology $H^{*}(G//H)=\oplus
H^{*}_{k}(G//H)$; this is the so-called {\it lower
grading} on $H^{*}(G//H)$.

According to~\cite[Theorem 2]{Hal},
(also cf.~\cite[Proposition 6.4]{Sin}),
$H^{*}_{k}(G//H)=0$ for $k>\mathrm{rank}(G)-\mathrm{rank}(H)$. In our case
$\mathrm{rank}(G)-\mathrm{rank}(H)=1$, and hence
$H^{*}_{k}(G//H)=0$ for $k>1$.
Next observe that
$H^{ev}_{1}(G//H)=0$ since $|q_{i}|$ is odd for any $i$ and
$H^{odd}(B_H)=0$. Similarly
$H^{odd}_{0}(G//H)=0$.
Therefore, $H^{ev}(G//H)=H^{*}_{0}(G//H)$ which is equal to the quotient of
$H^*(B_H)$ by $\bar{d}(\Lambda(q_1,...q_n))$. Since $H$ is a torus, $H^*(B_H)$
is generated by $2$-dimensional classes and hence by above the same
is true for $H^{ev}(G//H)$.

By Lemma~\ref{vanish} this implies that
$Der_{-}(H^{ev}(G//H),H^{*}(G//H))=0$ which by the splitting
criterion~\ref{prop: no-der-split-rigid}
means that $(G//H,T,k)$ is splitting rigid for any $T,k$.
\end{proof}

\section{Splitting rigidity for known positively curved manifolds}
\label{sec: pos curv}

\begin{prop}\label{positcurved}
Let $C$ be a known closed simply-connected positively curved
manifold, and $T$ be a torus.
If $\xi$ is a vector bundle over $C\times T$
such that $E(\xi)$ admits a complete metric with $\sec\ge 0$,
then $\xi$ virtually comes from $C$.
\end{prop}
\begin{proof}
First of all, note that all known even dimensional positively
curved manifolds belong to $\mathcal H$. Indeed,
it follows from the classification theorem of
positively curved homogeneous spaces~\cite{Wallach} that any even
dimensional homogeneous space belongs to
$\mathcal H(2n)\subset\mathcal H$ for some $n>0$.
The only known example of a positively curved even-dimensional
manifold which is
nondiffeomorphic to a positively curved homogeneous space
is the space $M^{6}=SU(3)//T^{2}$~\cite{Esc}. This space is an
$S^{2}$-bundle over $CP^{2}$ and, therefore, it lies in $\mathcal H$.

Now we establish splitting
rigidity for all known odd-dimensional positively curved
manifolds. Those are the standard spheres, the Berger
$7$-dimensional homology sphere $B^7=Sp(2)/Sp(1)_{\text{max}}$,
the Eschenburg  $7$-manifolds $E^7_{k,l,m,n}$~\cite {Esc1} all
obtained as biquotients of $SU(3)$ by $S^1$, and the Bazaikin
$13$-manifolds $B^{13}_{k,l,m,n}$~\cite{Baz} obtained as
biquotients of $SU(5)$ by $Sp(2)\times S^1$.

By~\ref{easy ex of split rigidity}, any odd-dimensional rational
homology sphere is splitting rigid, so it remains to deal with the
Eschenburg and Bazaikin manifolds.
By Lemma~\ref{esch-rational} below, all Eschenburg manifolds are
rationally homotopy equivalent to
$S^2\times S^5$, and all Bazaikin manifolds are rationally homotopy
equivalent to $CP^2\times S^9$. Since both $CP^2$ and $S^2$
belong to $\mathcal H (2)$, the proof of Theorem~\ref{generalproducts}
implies that any negative derivation of $H^*(S^2\times S^5)$ or
$H^*(CP^2\times S^9)$ vanishes on even cohomology.
Now Proposition~\ref{prop: no-der-split-rigid} implies
splitting rigidity of all Bazaikin and Eschenburg manifolds.
\end{proof}
\begin{lem}\label{esch-rational}
All Eschenburg manifolds are rationally homotopy equivalent to
$S^2\times S^5$; all Bazaikin manifolds are rationally homotopy
equivalent to $CP^2\times S^9$
\end{lem}
\begin{proof}
We only give a proof for the Bazaikin manifolds; the Eschenburg
manifolds are treated similarly. Let $B^{13}$ be a Bazaikin
manifold. From the homotopy sequence of the fibration
\[Sp(2)\times S^1\to SU(5)\to B^{13}\]
one easily sees that $B^{13}$ has the same rational homotopy
groups as $CP^2\times S^9$. Let $M=(\Lambda V,d)$ be the minimal
model of $B^{13}$. It is well-known (e.g. see~\cite[Theorem
15.11]{FHT}) that $V\cong \mathrm{Hom}_{\mathbb{Z}}(\pi_{*}(B),
\mathbb Q)$. Therefore, $M$ is a free graded algebra on generators
$x_2,y_5,y_9$ with the degrees of the generators given by the
subscripts. Obviously, $d(x_2)=0$ and $d(y_5)=kx_2^3$ for some
rational $k$. Note that $k\neq 0$, else we would have
$H^5(M)\cong\mathbb{Q}$ which is known not to be the case
by~\cite{Baz}. Replacing $y_5$ with $y_5/k$, we can assume that
$k=1$. It is clear that $dy_9$ must be equal to $lx_2^5$ for some
rational $l$. We claim that any such minimal model is isomorphic
to the one with $l=0$. Indeed, let $M_l=\langle
\Lambda(x_2,y_5,y_9)|dx_2=0,dy_5=x_2^3,dy_9=lx_2^5\rangle$.
Consider the map $M_0\to M_l$ given by $x_2\to x_2, y_5\to y_5,
y_9\to y_9-ly_5x_2^2$. This map is easily seen to be a
DGA-isomorphism with the inverse given by $x_2\to x_2, y_5\to y_5,
y_9\to y_9+ly_5x_2^2$. Since $M_0$ is a minimal model of
$CP^2\times S^9$, the proof is complete.
\end{proof}

\begin{rmk}
As we explained above, any known closed (simply-connected)
even-dimensional positively curved manifold $C$ belongs to
$\mathcal H(2n)$ for some $n$. The same is true for all
known~\cite{Wil2} even-dimensional manifolds with sectional
curvature positive on an open dense subset, such as projectivized
tangent bundles of $HP^n$ and $CP^n$. Therefore, all these
examples are also splitting rigid for any $T,k$.
\end{rmk}

\section{Splitting rigidity and derivations in minimal models}
\label{sec: der-in-min}
In this section, we study splitting rigidity
using methods of rational homotopy theory.
We prove that a triple $(C,T,k)$
is splitting rigid if, for any derivation of the minimal model
of $C$ that commutes with differential and has degree within
$[-\dim(T),0)$, the induced derivation on $H^*(C)$ vanishes
on $Char(C,k)$. The converse to this statement is proved
in~\ref{prop: split-rigid-no-der-in-model}
under various assumptions on $(C,T,k)$, such as
$2k\ge\dim(C\times T)+3$, or $p_i(TC)\in Char(C,k)$.

Thus, under either of the assumptions, the splitting rigidity is a
phenomenon of rational homotopy theory, in other words, whether or
not $(C,T,k)$ is splitting rigid depends only on $k$, $\dim(T)$,
and the minimal model (or equivalently, the rational homotopy
type) of $C$. This is no longer true for smaller $k$, in fact in
Section~\ref{sec: ex}, we give an example of two triples
$(C,T,6)$, $(M ,T,6)$ with  homotopy equivalent$C$ and $M$, such
that $(C,T,6)$ is splitting rigid, while $(M ,T,6)$ is not.

Let $(M_C,d_C)$, $(M_T,d_T)$ be (Sullivan) minimal models for $C$,
$T$. Since  $H^*(T)$ is a free exterior algebra, we can assume
that $M_T=H^*(T)$ and $d_T=0$. Then $(M_C\otimes M_T, d)$ is a
minimal model for $C\times T$, where for $x\in M_C$, $t\in M_T$
\[d(x\otimes t)=d_C(x)\otimes t+(-1)^{|x|} x\otimes d_T(t)=d_C(x)\otimes t.\]

First, we modify the arguments of Section~\ref{sec:
taylor-in-cohom} to produce the Taylor expansion in $M_C\otimes
M_T$ of any self-homotopy equivalence of $C\times T$. Let $h$ be a
self-homotopy equivalence of $C\times T$. Then the induced map of
minimal models $h^\#$ is an isomorphism~\cite[12.10(i)]{FHT}.

For $x\in M_C$, consider
$\frac{\partial h^\#}{\partial t_i}(x)\in M_C$
such that $h^\#(x\otimes 1)=\sum_{i}(1\otimes t_i)
(\frac{\partial h^\#}{\partial t_i}(x)\otimes 1)$.
Since $h^\#$ commutes with $d$ and $d(1\otimes t_i)=0$,
we get that each $\frac{\partial h^\#}{\partial t_i}$ commutes
with $d_C$, up to sign, and therefore induces a linear
self-map of $H^*(C)$.

As in Section~\ref{sec: taylor-in-cohom},
since $h^\#$ is an algebra isomorphism, the maps
$\frac{\partial h^\#}{\partial t_i}$
satisfy the same recursive identities,
obtained from $h^\#(xy)=h^\#(x)h^\#(y)$
by collecting the terms next to $1\otimes t_i$'s.
Thus, $\frac{\partial h}{\partial t_0}$ is an isomorphism of
$(M_C,d_C)$, and $\frac{\partial h^\#}{\partial t_0}\circ
(\frac{\partial h^\#}{\partial t_0})^{-1}$
is a derivation of $M_C$ of degree $-1$.
More generally, if $\frac{\partial h^\#}{\partial t_i}(x)=0$
for all $x\in M_C$ and all $0<i<k$,
then $\frac{\partial h^\#}{\partial t_k}\circ
(\frac{\partial h^\#}{\partial t_0})^{-1}$ is a derivation
of $M_C$ of degree $-|t_k|$.
As we noted before, $\frac{\partial h^\#}{\partial t_k}\circ
(\frac{\partial h^\#}{\partial t_0})^{-1}$ commutes with $d_C$, up to sign.
Hence, it induces a degree $-|t_k|$
derivation of $H^*(C)$.
A slight variation of the proof of
Proposition~\ref{prop: no-der-split-rigid},
implies the following.

\begin{prop}\label{prop: no-der-in-model-split-rigid}
If each negative  derivation of $H^*(C)$
induced by a derivation of $M_C$  of
degree $\ge -\dim(T)$
vanishes on $Char(C,k)$,
then $(C,T,k)$ is splitting rigid.
\end{prop}

\begin{rmk} The assumptions of
Proposition~\ref{prop: no-der-in-model-split-rigid}
only involve the minimal model of $C$.
Thus, if
$C^\prime$ is rationally homotopy equivalent to $C$ and
$(C,T,k)$ satisfies the assumptions of
Proposition~\ref{prop: no-der-in-model-split-rigid}, then
so does $(C^\prime, T, k)$.
\end{rmk}

The following Lemma shows that
an integer multiple of any negative
derivation of $M_C$ can be ``integrated'' to a self-homotopy
equivalence of $C\times T$. This gives many examples of triples
which are not splitting rigid, and is a crucial ingredient
in the proof of Theorem~\ref{mainex}.

\begin{lem} \label{prop: der-induce-he}
Let $D$ be a negative  derivation of $H^*(C)$
induced by a derivation of $M_C$ that commutes with $d_C$.
Let $T$ be a torus with $|D|\ge -\dim(T)$.
Then there are integers $i,m>0$, and a
self-homotopy equivalence $h$ of $C\times T$
such that $h^*(a\otimes 1)=a\otimes 1+(1\otimes t_i)(mD(a)\otimes 1)$
for $a\in H^*(C)$, and $h_{CC}\sim\mathrm{id}_C$,
$\pi_T\circ h=\pi_T$.
\end{lem}
\begin{proof}
Let $\tilde D$ be a derivation of $M_C$ that induces $D$.
The fact $\tilde D$ is a derivation implies that
the map $\phi\co M_C\to M_C\otimes M_T$ defined by
$\phi(x)=(x\otimes 1)+(1\otimes t_i)(\tilde D(x)\otimes 1)$
is a DGA-homomorphism.

Being a DGA-homomorphism,
$\phi$ defines a (unique up to homotopy) map of rationalizations
$f\co C_0\times T_0\to C_0$ where we can
choose $f$ so that $f\circ i_{C_0}=\mathrm{id}_{C_0}$.
Look at the diagram
\begin{equation*}
\xymatrix{
C\times T \ar@{-->}[d]^{p} \ar@{-->}[dr]^{\tilde{f}} &\\
C\times T\ar[d]^{r}  & C\ar[d]^{r_C}\\
C_0\times T_0\ar[r]^{f}&C_0}
\end{equation*}
where $r=r_T\times r_C\co C\times T\to C_0\times T_0$ is the
rationalization, and try to find a finite covering $p$, and a map
$\tilde f$ that makes the diagram commute, and satisfies $\tilde
f\circ i_C=\mathrm{id}_C$.

It follows from an obstruction theory argument as 
in~\cite[Section 4]{BK1}  
that such $p$, $\tilde f$ can be constructed, where the
key point is that all obstructions are torsion since the homotopy
fiber of $r$ has torsion homotopy groups, and all torsion
obstruction vanish after precomposing with a suitable finite cover
$p$. Furthermore (cf. the proof of Lemma~\ref{lem: bundles}), one can
choose $p,\tilde{f}$ satisfying $\tilde{f}\circ
i_C=\mathrm{id}_C$, $p=\mathrm{id}_C\times (\times^n)$, where
$\times^n\co T\to T$ is the $n$-th power map, so that, for some
positive integer $m$, the induced map on $H^*(C)$ satisfies
\[
\tilde f^*[a\otimes 1]=[a\otimes 1]+m(1\otimes t_i)([\tilde
D(a)]\otimes 1).
\] 
Finally, by Whitehead's theorem, the map
$h(c,t)=(\tilde f(c,t),t)$ is a homotopy equivalence with the
desired properties.
\end{proof}

\begin{prop}\label{prop: split-rigid-no-der-in-model}
Assume that $(C,T,k)$ is splitting rigid and either
$2k\ge\dim(C\times T)+3$, or $p_i(TC)\in Char(C,k)$ for all $i>0$.
If $D$ is a negative  derivation of $H^*(C)$ of degree $\ge -\dim(T)$
induced by a derivation of $M_C$, then
$D$ vanishes on $Char(C,k)$.
\end{prop}
\begin{proof}
Arguing by contradiction, let $D$ be a negative derivation
of $H^*(C)$ with $|D|\ge -\dim(T)$
such that $D$ does not vanish on $Char(C,k)$.

First, we show that $D$ is nonzero on either Euler or Pontrjagin
class of a rank $k$ bundle $\xi_C$ over $C$. Assume, say that
$k=2m$ so that $Char(C,k)=(\oplus_{i=1}^{m-1}H^{4i}(C))\oplus
H^{2m}(C)$. If $a\in H^{2m}(C)$ satisfies $D(a)\neq 0$, then by
Lemma~\ref{lem: bundles}, $a$ is proportional to the Euler class
of some rank $k$ vector bundle $\xi_C$ over $C$. If $a\in
H^{4i}(C)$, then by Lemma~\ref{lem: bundles}, $a$ is proportional
to the $i$th Pontrjagin class of some bundle of rank $k$ vector
bundle $\xi_C$ over $C$.

Say, suppose that $D$ is nonzero on the Euler class $e(\xi_C)$.
By Lemma~\ref{prop: der-induce-he}, there is a
self-homotopy equivalence $f$ of $C\times T$ such that
$f^*(e(\xi_C)\otimes 1)$ does not lie in $H^*(C)\otimes 1$.
Now there are two cases to consider.

If $2k\ge\dim(C\times T)+3$, then by Haefliger's embedding
theorem~\cite{Hae} the homotopy equivalence $f\co C\times T\to
E(\xi_C)\times T$ is homotopic to a smooth embedding $q$. The
normal bundle $\nu_q$ has rank $k$. Since $q$ and $f$ are
homotopic, $\nu_q$ and $q^\#(\xi_C\times T)$ have equal Euler and
Pontrjagin classes. This follows from the intersection pairing
interpretation of the Euler class, and the Whitney sum formula for
the total Pontrjagin class (see~\cite[Section 3]{BK1} for
details). Thus, $e(\nu_q)=f^*(e(\xi_C)\otimes 1)$ does not lie in
$H^*(C)\otimes 1$. So $\nu_q$ does not virtually come from $C$,
while $E(\nu_q)$ is diffeomorphic to $T\times E(\xi_C)$, and this
means that $(C,T,k)$ is not splitting rigid.

If $p_i(TC)\in Char(C,k)$ for all $i>0$, then $D(p(TC))=0$. In
other words, $f^*p(TC)=p(TC)$, so by Lemma~\ref{lem: surgery},
there is a diffeomorphism $q$ of $C\times T$ such that
$q^*(e(\xi_C)\otimes 1)$ does not lie in $H^*(C)\otimes 1$. Then
the pullback bundle $f^\#(\xi_C\times T)$ does not virtually come
from $C$, while its total space is diffeomorphic to $T\times
E(\xi_C)$; thus $(C,T,k)$ is not splitting rigid.
\end{proof}
\begin{rmk}
Note that if $k\ge \dim(C)$, then $p_i(TC)\in Char(C,k)$ for all
$i>0$.
\end{rmk}
\begin{rmk}
The proof of Proposition~\ref{prop: split-rigid-no-der-in-model}
implies the weak converse of the Proposition~\ref{characterization}:
if $(C,T,k)$ is splitting rigid, and either
$2k\ge\dim(C\times T)+3$, or $p_i(TC)\in Char(C,k)$,
then any homotopy equivalence of $C\times T$ maps
$Char(C,k)\otimes 1$ to itself.
\end{rmk}
The following example is due to T.~Yamaguchi~\cite[Example 3]{TYam}
and it was constructed in response to a question in the previous
version of this paper.

\begin{ex}\label{ex:Yam}
Consider a minimal Sullivan algebra $M$  with the generators
$\{x, y, z, a, b, c \}$, whose degrees are respectively
$2, 3, 3, 4, 5, 7$, and the differential given by $d(x)=d(y)=0$,
$d(z)=x^2$, $d(a)=xy$, $d(b)=xa+yz$, $d(c)=a^2+2yb$.
A direct computation shows that the rank of $H^i(M)$ is equal to
$1$ for $i =0,2,3,11,12,14$, is equal to $2$ for $i =7$, and is equal
to 0 for other $i$'s. Also the generators of $H^*(C)$ as an algebra
are \[ x,\ y,\ e=ya,\ f=xb-za,\ g=x^2c-xab+yzb,\ h=3xyc+a^3. \] The
products are all trivial, except $xh$, $yg$, and $ef$ in $H^{14}(M)$.
Thus $H^*(M)$ satisfies the Poincare duality, and
therefore by rational surgery~\cite[Theorem 13.2]{Sul},
there exists a closed simply
connected elliptic $14$-dimensional manifold
$C$ with minimal model $M$.

Then $\Der (H^*(C))$ is not zero since there is, for example, a
non-zero derivation  $(g,y)$ of degree $-8$. Here $(p,q)$ stands
for the derivation which send p to q and other generators to zero.
One can check that $(g,y)$ is indeed a derivation. On the other
hand, another direct computation (see~\cite[Example 3]{TYam})
shows that all derivations induced from $M$ vanish on $H^*(C)$.
Therefore, the manifold $X=C\times C\times C\times C$ has the same
property, but there exists a derivation $D$ of $H^*(X)$ of degree
$-8$ which is non-zero on $H^{44}(X)$. Thus, $(X, T, k)$ is
splitting rigid for any $T,k$, but this fact can not  be seen by
looking only at $H^*(X)$.
\end{ex}
The following lemma, combined
with Proposition~\ref{prop: no-der-in-model-split-rigid},
shows that the property of being
``splitting rigid for all $k$'' depends only on $\dim(T)$
and the rational homotopy type of $C$.
\begin{lem}
Let $(C,T,k)$ be splitting rigid for all $k$.
Then if $D$ is a negative  derivation of $H^*(C)$ of degree $\ge -\dim(T)$
induced by a derivation of $M_C$, then
$D$ vanishes on $H^{even}(C)$.
\end{lem}
\begin{proof}[Sketch of the proof]
As in the proof of
Proposition~\ref{prop: split-rigid-no-der-in-model}, use $D$
to construct a self-homotopy equivalence $f$ of $C\times T$.

If $D$ is nonzero on $H^{4i}(C)$ for some $i$, then $D$ is nonzero
on the $i$th Pontrjagin class of some bundle $\xi_C$ over $C$.
Then for some large $k$, the map $f\co C\times T\to E(\xi_C)\times T$
is homotopic to a smooth embedding.
By the same argument as in the proof of
Proposition~\ref{prop: split-rigid-no-der-in-model}, the normal bundle
to this embedding does not virtually come from $C$ and, therefore,
$(C,T,k)$ is not splitting rigid.

If $D$ vanishes on $\oplus_{i}H^{4i}(C)$, then $D(p(TC))=0$
so $f$ preserves $p(TC)$. Hence, by Lemma~\ref{lem: surgery},
replacing $f$ with some power of $f$,
we can assume  that $f$ is homotopic to
a diffeomorphism. If $D$ is nonzero on $H^{2i}(C)$ for some $i$,
then $D$ is nonzero on the Euler class of some bundle $\xi_C$ over $C$.
Looking at the bundle $f^\#(\xi_C\times T)$ shows that
$(C,T,2i)$ is not splitting rigid.
\end{proof}

\section{Proof of Theorem~\ref{mainex}}
\label{sec: non-splil-rigid-ex}
Let $G=SU(6)$ and $H=SU(3)\times SU(3)$.
According to~\cite[Chapter5, Example 4.14]{OT}
(cf.~\cite[Section 11.14]{GHV} and~\cite[Proposition 5.16]{FHT}),
the minimal model of
$G/H$ is given by $(M, d)=(\Lambda(y_4,y_6, x_7, x_9, x_{11}), d)$
with the degrees given by the subscripts, and $d(y_4)=0=d(y_6)$,
$d(x_7)=y_4^2$, $d(x_9)=2y_4y_6$, $d(x_{11})=y_6^2$.

Now it is straightforward to compute the cohomology algebra of $G/H$.
In particular, $G/H$ has nonzero Betti numbers
only in dimensions $0, 4, 6, 13, 15, 19$
and the cohomology groups in dimensions $4, 6$ are
generated by the classes $[y_4], [y_6]$.

Let $\xi$ be the rank $3$ complex vector bundle over $G/H$
classified by $p\co G/H\to BSU(3)$ which is the composition of the
classifying map $G/H\to BH=BSU(3)\times BSU(3)$ for the bundle
$G\to G/H$ with the projection $BSU(3)\times BSU(3)\to BSU(3)$ on
the first factor. Since $G/(SU(3)\times 1)$ is $6$-connected, from
the Serre spectral sequence of the bundle $ G/SU(3)\to G/H\to
BSU(3)$ we see that the map $p^{*}\co H^{i}(BSU(3))\to H^{i}(G/H)$
is an isomorphism for $i\le 6$. Hence, $c_3(\xi)$ is the generator
of $H^{6}(G/H)$, and by rescaling $y_6$, we can assume that
$c_3(\xi)=[y_6]$. Note that $c_3(\xi)$ is equal to the Euler class
$e(\xi_\mathbb R)$ of $\xi_\mathbb R$, the realification of $\xi$.

Alternatively, $\xi$ can be described as the associated bundle
to the principal bundle $G\to G/H$
via the representation $\rho$ of $SU(3)\times SU(3)$ given by
the projection onto the
first factor followed by the standard action of $SU(3)$ on
$\mathbb C^{6}$.  Thus, $E(\xi)$ admits a complete metric with
$\sec\ge 0$ such that the zero section is a soul.

To finish the proof it remains to find a torus $T$, and a
self-diffeomorphism $f$ of $C\times T$ such that
$f^*(c_3(\xi)\otimes 1)\notin H^*(G/H)\otimes 1$, because
then the bundle $f^\#\xi_\mathbb R$ does not virtually come
from $C$, and $E(f^\#\xi_\mathbb R)$ carries a complete metric
with $\sec\ge 0$ and zero section being a soul.

Because $M$ is free, the linear map $\tilde D\co M\to M$ defined
by $\tilde D(y_4)=\tilde D(x_7)=0$, $\tilde D(y_6)=y_4$, $\tilde
D(x_9)=x_7$, $\tilde D(x_{11})=2x_9$ is a derivation of $M$ of
degree $-2$ (see~\cite[page 141]{FHT}). By computing on the
generators, it is straightforward to see that $\tilde D$ commutes
with $d$, and hence induces a derivation $D$ of $H^*(G/H)$ such
that $D([y_6])=[y_4]$.

By Lemma~\ref{prop: der-induce-he}, there is a positive integer $m$,
and a self-homotopy equivalence $h$ of $C\times T$ where $\dim(T)=2$
such that $\pi_T\circ h=\pi_T$ and
$h^*(a\otimes 1)=a\otimes 1+(1\otimes t_3)(mD(a)\otimes 1)$
for any $a\in H^*(G/H)$.

Note that $h^*$ preserves the total Pontrjagin class
of the tangent bundle to $G/H\times T$. Indeed, since $T$
is parallelizable and $H^{4i}(G/H)$ are only nonzero if $i=0,1$,
it suffices to show that $h^*$ preserves $p_1(G/H)\otimes 1$, or
equivalently, that $mD$ vanishes on $p_1(G/H)$. In fact, more is true,
namely, $D$ vanishes on $H^4(G/H)$ since $D([y_4])=[\tilde D(y_4)]=0$.
By Lemma~\ref{lem: surgery} below, some power
\[f=h^k=\underbrace{h\circ\dots\circ h}_k\] of $h$ is homotopic to a
diffeomorphism. Since $\pi_T\circ h=\pi_T$, we have $\pi_T\circ f=\pi_T$
so that $f^*(1\otimes t)=1\otimes t$ for any $t\in H^*(T)$.
Combining with $(t_3)^2=0$, we get for $a\in H^*(C)$
\[f^{*}(a\otimes 1)=a\otimes 1+(1\otimes t_3)(kmD(a)\otimes 1).\]
Since $D(c_3(\xi))=D([y_6])=[y_4]\neq 0$, we get
$f^*(c_3(\xi)\otimes 1)\notin H^*(G/H)\otimes 1$,
and the proof of Theorem~\ref{mainex} is complete.

\section{Proving splitting rigidity without rational homotopy}
\label{sec: ex}
This section contains an example of two triples $(C,T,6)$,
$(M,T,6)$ such that $C$ and $M$ are homotopy equivalent and
$(C,T,6)$ is splitting rigid, while $(M, T,6)$ is not.

Let $M=S^3\times S^3\times S^{10}\times S^{11}$.
Note that $(M,T,6)$ is not splitting rigid if $\dim(T)\ge 3$.
Indeed, let $\xi_M$ be the pullback of $TS^6$ via the map
$M\to S^6$, which is the composition of the projection
$M\to S^3\times S^3$ followed by a degree one map
$S^3\times S^3\to S^6$. Then $\xi_M$ has nonzero Euler class.
If $\dim(T)\ge 3$, then one can easily construct a
self-diffeomorphism $f$ of $M\times T$ such that
$f^* e(\xi_M\times T)$ does not lie in $H^6(M)\otimes 1$,
in particular, the bundle $f^\#(\xi_M\times T)$ does not
virtually come from $M$.

\begin{prop}\label{prop: ld-ex}
There exists a smooth manifold $C$ homotopy equivalent
to $M$ such that $(C,T,6)$ is splitting rigid
for any $T$.
\end{prop}
\begin{proof}
By a surgery argument as in~\cite[A.1]{BK0}, there is a closed
smooth manifold $C$ which is homotopy equivalent to $M$
and has $p(TC)=1+p_4(TC)$ with $p_4(TM)\neq 0$.

To check that $(C,T,6)$ is splitting rigid, we
need to start with an arbitrary rank $6$ vector bundle $\xi$ over $B=C\times T$
that satisfies ($\ast$), and prove that
$\xi$ virtually comes from $C$, or equivalently, that
the Euler and Pontrjagin classes of $\xi$ lie in
$Char(C,6)\otimes 1$.

We shall  borrow notations and arguments from the proof
of Proposition~\ref{characterization}. As in~\ref{characterization},
we can assume that $E(\xi)$
is the total space of a vector bundle $\eta$ which is the product of
$T$ and a vector bundle $\eta_{C^\prime}$
over a closed smooth simply-connected manifold $C^\prime$.
Let $S=C^\prime\times T$ and $g\co B\to S$ by the homotopy equivalence
as in~\ref{characterization}.

Note that $g^*p(TS)=p(TB)$. Indeed, as in~\ref{characterization},
$g^*p(TN|_S)=p(TN|_B)$. So $g^*p(\eta)g^*p(TS)=p(\xi)p(TB)$.
The only Pontrjagin classes of
a rank $6$ vector bundle that have
a chance of being nonzero are $p_1$, $p_2$, $p_3$.
Since $C^\prime$ has zero cohomology in dimensions
$4$, $8$, $12$, we get $p(\eta)=1$.
By the same argument, $p(TS)=1+p_4(TS)+p_6(TS)$.
We have
\[
g^*(1+p_4(TS)+p_6(TS))=(1+p_1(\xi)+p_2(\xi)+p_3(\xi))(1+p_4(TB)),
\]
hence $p_6(TS)=0=p_i(\xi)$ for all $i$, and $g^*p(TS)=p(TB)$.
Then one easily sees that $g_{CC^\prime}$ maps $p(TC^\prime)$
to $p(TC)$. Hence, $g_{CC^\prime}\times g_{TT}$ maps
$p(TS)$ to $p(TB)$, and therefore,
$h=(g^{-1}_{CC^\prime}\times g^{-1}_{TT})\circ g$
preserves $p(TB)$.

As in~\ref{characterization},
it suffices to show that $h$ preserves
$Char(C,6)\otimes 1= H^6(C)\otimes 1$.
We think of $H^*(C)$ as an exterior algebra on generators
$x,y, q, s$
corresponding to spheres $S^3$, $S^3$, $S^{10}$, $S^{11}$
so we need to show that $h^*(xy)=xy$.
By rescaling we can assume that $p(TB)=xyq\otimes 1$.

By dimension reasons
$\frac{\partial h^*}{\partial t_i}(xy)=0$ unless
$|t_i|$ is $3$ or $6$. Similarly,
$\frac{\partial h^*}{\partial t_i}(q)=0$
unless $|t_i|$ is $4$, $7$ or $10$.
Now collecting terms next to
$1\otimes t_i$'s in the identity $xyq\otimes 1=h^*(xyq\otimes 1)$,
we conclude that $q\frac{\partial h^*}{\partial t_i}(xy)=0$,
and hence $\frac{\partial h^*}{\partial t_i}(xy)=0$
for all $i>0$. Thus, $h^*(xy)=xy$ as promised.
\end{proof}

\begin{rmk}
The above example shows that
Proposition~\ref{prop: split-rigid-no-der-in-model} fails
without assuming either
$2k\ge\dim(C\times T)+3$, or $p_i(TC)\in Char(C,k)$.
Indeed, $(H^*(C), 0)$ is a minimal model of $C$
and there is a degree $-3$ derivation of $H^*(C)$ given by
$D(x)=1$, $D(y)=D(q)=D(s)=0$ which does not vanish
on $Char(C,6)=H^6(C)$. Namely, $D(xy)=y$.
Yet $(C,T,6)$ is splitting rigid.

It is instructive to see where the proof of
Proposition~\ref{prop: split-rigid-no-der-in-model} fails.
Using $D$, we produce a self-homotopy equivalence
$f$ of $C\times T$, and an $\mathbb R^6$-bundle $\xi_C$
with $xy=e(\xi_C)\otimes 1$. However, the homotopy
equivalence $f\co C\times T\to E(\xi_C)\times T$
is not homotopic to a smooth embedding.
\end{rmk}

\begin{rmk}
As always with splitting rigid triples,
the total spaces of ``most'' $\mathbb R^6$-bundles
over $C\times T$ do not admit complete metrics with $\sec\ge 0$.
We do not know whether $C$ in Proposition~\ref{prop: ld-ex}
admits a metric with $\sec\ge 0$. Yet, no currently
known method rules out the existence of $\sec\ge 0$ on $C$,
because $C$ is homotopy equivalent to a closed nonnegatively
curved manifold, and $C$ admits a metric of positive scalar
curvature, for $C$ is spin and $\dim(C)=27\equiv 3\ (\mathrm{mod}\ 8)$,
so~\cite{Sto} applies.
\end{rmk}

\section{Nonnegatively curved vector bundles with souls equal to the
zero sections}
\label{sec: souls-sections}
The purpose of this section is to obtain restrictions on normal bundles
to souls in nonnegatively curved manifolds.
In other words, we look for conditions on a vector bundle $\xi$
ensuring that $E(\xi)$ admits no complete nonnegatively curved metric
such that the zero section is a soul.
The assumption that a given submanifold is a soul
imposes a nontrivial restriction on the metric, so it is no surprise that
we get stronger results on obstructions.

Our exposition is parallel to the one in Section~\ref{sec: split criterion}.
We say that a vector bundle $\xi$ over $C\times T$  satisfies
condition ($\ast\ast$) if

\hfill $\begin{array}{l} \text{there is a finite cover } \pi\co
C\times T\to C\times T,\ \text{a closed manifold}\ C^\prime,\
\text{and}\\ \text{a diffeomorphism }f\co C^\prime\times T\to C\times T
\text{ such that the bundle }f^\#\pi^{\#}(\xi)\\ \text{virtually
comes from}\ C^\prime.
\end{array}$\hfill $(\ast\ast)$
\begin{ex} According to Theorem~\ref{geomsplit},
if $S$ is a soul in a complete
nonnegatively curved manifold, then the normal bundle to $S$ satisfies
($\ast\ast$).
\end{ex}
{\bf Caution.} Clearly, if $\xi$ satisfies ($\ast\ast$),
it also satisfies condition ($\ast$) from
Section~\ref{sec: split criterion}.
The converse is generally false, as
the following example shows.
\begin{ex}
Let $M=S^{3}\times S^{5}\times S^{7}$. By a surgery argument as
in~\cite[A.1]{BK0}, there is a closed $15$-dimensional manifold $C$
which is homotopy equivalent to $M$ and such that $p_{2}(TC)$,
$p_{3}(TC)$ are nonzero.
Clearly, there is a derivation of $H^{*}(C)$ of degree $-3$
which is nonzero on $H^{8}(C)$. Since $(H^{*}(C),0)$ is the minimal
model for $C$, Proposition~\ref{prop: split-rigid-no-der-in-model}
implies that $(C,T,15)$ is not splitting rigid if $\dim(T)\ge 3$.
Thus, there is a rank $15$ bundle $\xi$ over $C\times T$
which satisfies ($\ast$) but which does not virtually come from $C$.
It remains to show that $\xi$ does not satisfy ($\ast\ast$).
If it does, then after passing to a finite cover, we can assume that,
for some diffeomorphism $f\co  C^\prime\times T\to C\times T$,
$f^\#\xi$ virtually comes from $C^{\prime}$.
Since $f$ is a diffeomorphism and $T$ is parallelizable,
$f^*(p_i(TC)\otimes 1)=p_i(TC^\prime)\otimes 1$ for all $i$.
Since $p_i(TC)$ generates $H^{4i}(C)$ and $f^{*}$ is an isomorphism,
this means that
$f^{*}(H^{4i}(C)\otimes 1) = H^{4i}(C^{\prime})\otimes 1$
for all $i$. Since $f^*p_i(\xi)=p_{i}(f^{\#}(\xi))\in
H^{4i}(C^{\prime})\otimes 1$
for all $i$, this implies that $p_i(\xi)\in H^{4i}(C)\otimes 1$
for all $i$
and hence $\xi$ virtually comes from $C^{\prime}$. This is a
contradiction, so $\xi$ cannot satisfy ($\ast\ast$).
\end{ex}
Recall that, given a compact Lie group $G$,
a rank $n$ vector bundle $\xi$ has structure group $G$ if $\xi$
is associated with a principal $G$-bundle via some representation
$G\to O(n)$.
We now have the following splitting criterion similar to
Proposition~\ref{characterization}:
\begin{prop}\label{smoothchar}
Let $\xi$ be a vector bundle over $C\times T$ which has
structure group $O(k)$ and satisfies ($\ast\ast$).
If any self homotopy equivalence of $C\times T$ maps
$Char(C,k)\otimes 1$ to itself,
then $\xi$ virtually comes from $C$.
\end{prop}
\begin{proof} Passing to a finite cover,
we can assume that $\xi$ is orientable, and that $\eta=f^{\#}\xi$
virtually comes from $C^{\prime}$,
where $f\co C^\prime\times T\to C\times T$ is a diffeomorphism.
As before, to show that $\xi$ virtually comes from
$C$ it is enough to check that its rational characteristic classes
lie in $H^{*}(C)\otimes 1$.
Proceeding exactly as in the proof of
Proposition~\ref{characterization}, we conclude that
$(f^{*})^{-1}$ maps $Char(C^{\prime},k)\otimes 1$ to $Char(C ,k)\otimes 1$.

Next note that all $e(\eta)$, $p_{i}(\eta)$ lie in the
subalgebra generated by $Char(C^{\prime},k)\otimes 1$.
Indeed, $\eta$ has the structure group $SO(k)$, so
$\eta$ is a pullback of a bundle over $BSO(k)$.
Since the cohomology of $BSO(k)$ is generated by $Char(BSO(k),k)$,
by naturality of the characteristic classes,
we see that $e(\eta), p_{i}(\eta)\in \langle Char(C^{\prime},k)\rangle$
for any $i$.

Since $(f^{*})^{-1}(\langle Char(C^{\prime},k)\otimes 1\rangle)\subset
\langle Char(C,k)\otimes 1\rangle$,
we conclude that all the characteristic classes of $\xi$ lie
in $\langle Char(C,k)\otimes 1\rangle\subset H^{*}(C)\otimes 1$,
hence $\xi$ virtually comes from $C$.
\end{proof}

Now all the splitting rigidity results that relied on
Proposition~\ref{characterization}, can be adapted to
this new setting. In particular, we obtain
\begin{thm}\label{thm: rk2-smooth-split-rigid}
Let $\xi$ be a vector bundle over $C\times T$ with
structure group $O(2)$. If
$E(\xi)$ admits a complete metric with $\sec\ge 0$ such that the zero
section is a soul, then $\xi$ virtually comes from $C$.
\end{thm}

\begin{thm}\label{smooth-biquotients}
Let $C=G//H$ be a simply connected biquotient of
compact  Lie groups such that $H$ is
semi-simple.
Let $\xi$ be a vector bundle over $C\times T$ whose structure group
can be reduced to a subgroup of $O(4)$.
If $E(\xi)$ admits a complete metric with $\sec\ge 0$ such that the zero
section is a soul,
then $\xi$ virtually comes from $C$.
\end{thm}

\section{Open problems}
\label{sec: prob}
\subsection{Induced derivations}
As we proved in Section~\ref{sec: der-in-min}, for sufficiently
large $k$, splitting rigidity of $(C,T,k)$ is equivalent to
vanishing on $Char(C,k)$ of all  negative derivations of degree
$\ge -\mathrm{dim}(T)$ induced from the minimal model $M_{C}$. Yet
all our geometric applications are proved by checking the stronger
condition that {\it all} negative derivations vanish on
$Char(C,k)$. This is mostly due to the fact that we do not know
how to effectively check which derivations of $H^{*}(C)$ are
induced from the minimal model.

Let us restate the problem in  purely rational homotopy theoretic
terms. It is well known that the space $\mathrm{Der}(M_{C})$ is a
differential graded Lie algebra (DGLA) with the differential given
by $D_{C}=[-,d_{C}]$. It is trivial to check that closed
derivations preserve $\mathrm{ker}(d_{C})$ and exact ones send
$\mathrm{ker}(d_{C})$ to $\mathrm{Im}(d_{C})$ (see~\cite{Gri} for
details). Therefore, we have a natural graded Lie algebra
homomorphism $m\co H^{*}(\mathrm{Der}(M_{C}))\to
\mathrm{Der}(H^{*}(C))$. We seek to understand the image of the
map $m$. Example~\ref{ex:Yam} produces an elliptic smooth manifold
$C$ such that $\mathrm{Im}(m)=0$ but
$\mathrm{Der}(H^{*}(C))\ne 0$.

To relate to our geometric applications we would like to find such
examples when $C$ is nonnegatively curved. It is easy to see that
$m$ is onto if $C$ is formal, but that is all we can generally say
at the moment. (Recall that a space $X$ is called formal if $X$ and
$(H^*(X),d=0)$ have isomorphic minimal models).

Let us also mention that according to Sullivan~\cite{Sul}, the DGLA
$(\mathrm{Der}(M_{C}),D_{C})$ is a
(Quillen) Lie algebra model for $Baut_{1}(C)$
(the classifying space for the identity component of the monoid of
self-homotopy equivalences of $C$),
and therefore, understanding the map $m$ can be helpful for computing the
rational homotopy groups of $Baut_{1}(C)$.

\subsection{Halperin's conjecture for biquotients}
As we mentioned in the introduction, the conjecture of Halperin that
any elliptic space of positive Euler characteristic belongs to
$\mathcal H$ has been verified for all homogeneous spaces of compact Lie
groups~\cite{ST}. However, it remains open for the natural bigger
class of elliptic spaces formed by biquotients.
According to~\cite{Sin}, a biquotient $G//H$ has a positive Euler
characteristic iff $\mathrm{rank} (G)= \mathrm{rank} (H)$.
Therefore, we pose the following
\begin{prob} Prove that any biquotient  $G//H$
belongs to $\mathcal H$ if $\mathrm{rank} (G)= \mathrm{rank} (H)$.
\end{prob}
This is unknown even for the simplest examples such
as $Sp(1)\backslash Sp(n)/SU(n)$.
\subsection{Nonnegatively curved vector bundles over rational H-spaces}

Most explicit examples of nonnegatively curved bundles are given by
{\it homogeneous} vector bundles, i.e. by vector bundles
associated to principal $H$-bundles $H\to G\to G/H$ via
some representations $H\to O(k)$.
For any given $k$, there are only
finitely many rank $k$ homogeneous vector bundles
(because the number of nonequivalent
irreducible representations $H\to O(k)$ is finite).
However, homogeneous vector bundles can fill a substantial
part of $[G/H,BO]$, the set of stable equivalence classes
of bundles over $G/H$.

If $H=1$, or more generally if
$G/H$ is rationally homotopy equivalent to
the product of odd-dimensional spheres,
then $H^*(BH)\to H^*(G/H)$ has trivial image~\cite[page 466]{GHV},
and hence, any homogeneous vector bundle over $G/H$ has zero
Euler and Pontrjagin classes.
Motivated by the above discussion, we pose the following:
\begin{prob}\label{newbundles} Does there exist a nonnegatively
curved vector bundle $\xi$ over a closed manifold
$C$ such that $C$ is rationally homotopy equivalent to the product
of odd-dimensional spheres, $\sec(C)\ge 0$, and
$e(\xi)\neq 0$ or $p_{i}(\xi)\ne 0$ for some $i>0$?
\end{prob}

\subsection {Nonnegatively curved nonsplitting rigid examples}
The example of a nonnegatively curved vector bundle $\xi$ that
satisfies condition ($\ast$) but does not virtually come from $C$
provided by Theorem~\ref{mainex} is essentially the only example
of this kind known to us. This is not very satisfactory, say,
because this example is unstable; that is, $\xi\oplus\epsilon^{1}$
does virtually come from $C$. More importantly, we want to
understand how ``generic'' such examples are.

To construct a stable example,
it suffices to find a vector bundle $\xi$ over $C$
with $\sec(E(\xi ))\ge 0$ and soul equal to the zero section,
and a negative derivation $D$ of $H^{*}(C)$ induced by a
derivation of the minimal model
such that $D(p(\xi))\ne 0$ and $D(p(TC))=0$.
While we think that many such examples exist, finding
an explicit one, say among homogeneous vector bundles,
seems to be an unpleasant task because
\begin{enumerate}
\item Pontrjagin classes of homogeneous vector bundles are often
difficult to compute, and
\item there is no easy algorithm for computing the space of negative
derivations of $H^{*}(G/H)$ or of the minimal model of $G/H$.
\end{enumerate}
One of the few cases when $\mathrm{Der}_{-}(H^{*}(G/H))$ is easily
computable is when $G/H$ is formal. According to~\cite[Theorem 12.2]{Oni2},
any formal compact homogeneous space $G/H$
is rationally homotopy equivalent
to the product of odd-dimensional spheres and an elliptic space $X$ of
positive Euler characteristic. Again by~\cite[Theorem 12.2]{Oni2},
the image of the homomorphism $H^{*}(BH)\to H^{*}(G/H)$ is equal to
the $H^{*}(X)$-factor.
If $X\in\mathcal H$ (i.e. if Halperin's conjecture holds for $X$),
then by Lemma~\ref{decomp}, one concludes
that if $\xi$ is a  homogeneous vector bundle and $D$ is a
negative derivation of $H^*(G/H)$, then $D$ vanishes
on $e(\xi)$, $p_i(\xi)$ for $i>0$.
Thus, if Halperin's conjecture is true, then
homogeneous vector bundles over formal homogeneous
spaces cannot be used to prove an analog of
Theorem~\ref{mainex}.

Finally, note that a positive solution to
Problem~\ref{newbundles} (for $C$ with $p(TC)=1$ which includes
the case when $C$ is a compact Lie group or the product of
odd-dimensional spheres)
yields an analog of Theorem~\ref{mainex},
because for any nontrivial element $a$ of $H^*(C)$,
there exists a negative derivation $D$ of $H^*(C)$ with
$D(a)\neq 0$, and $D(p(TC))$ vanishes by assumption.

\appendix
\section{Surgery-theoretic lemma}
We are grateful to Ian Hambleton for sketching the proof of the
following lemma.
\begin{lem}\label{lem: surgery}
Let $C$ be a closed smooth simply-connected manifold, and $T$
be a torus such that $\dim(C)+\dim(T)\ge 5$.
If $h$ is a self-homotopy equivalence of $C\times T$ that
preserves the rational total Pontrjagin class of $C\times T$, then
$h^m$ is homotopic to a diffeomorphism for some $m>0$.
\end{lem}
\begin{proof} Let $k=\dim(T)$, $n=\dim(C)$,
$I=[0,1]$, $B=C\times T$, $\pi=\pi_1(B)$
so that $\pi\cong\mathbb Z^k$.
Look at the following commutative diagram whose rows are
the smooth and the topological surgery exact sequences:
\[
\xymatrix{
[B\times I\ \mbox{rel}\ \partial,G/O]\ar[r]\ar[d]&
L_{n+k+1}(\mathbb Z\pi)\ar[r]\ar@{=}[d]& S_O(B)\ar[r]\ar[d]&
[B, G/O]\ar[r]\ar[d]& L_{n+k}(\mathbb Z\pi)\ar@{=}[d]\\
[B\times I\ \mbox{rel}\ \partial,G/Top]\ar[r]&
L_{n+k+1}(\mathbb Z\pi)\ar[r] & S_{Top}(B)\ar[r]&
[B, G/Top]\ar[r]& L_{n+k}(\mathbb Z\pi)}
\]
First, note that
$[B\times I\ \mbox{rel}\ \partial,G/Top]\to L_{n+k+1}(\mathbb Z\pi)$ is onto.
Indeed, by the Poincare duality with $L$-theory coefficients
$[B\times I\ \mbox{rel}\ \partial,G/Top]=H^0(B\times I;\mathbf{L})\cong
H_{n+k+1}(B\times I;\mathbf{L})$, so it suffices to show that
the homology assembly map (see e.g.~\cite[p216]{Dav})
\[ A_{n+k+1}\co H_{n+k+1}(B\times I;\mathbf{L}) \to
L_{n+k+1}(\mathbb Z\pi)\] is onto.
By naturality of the assembly and since $T$ is the classifying space
for $\pi_1(B\times I)$, $A_{n+k+1}$ factors as the composition of
the map $H_{n+k+1}(B\times I;\mathbf{L})\to H_{n+k+1}(T ;\mathbf{L})$
induced by the projection $B\times I\to T$, and the universal assembly
$H_{n+k+1}(T ; \mathbf{L})\to L_{n+k+1}(\mathbb Z\pi)$.
The former map is onto, since it has a section induced by a section of
$B\times I\to T$, while the latter map is an
isomorphism since $\pi\cong\mathbb Z^k$~\cite[Chapter 15B]{Wal}.
Thus, $A_{n+k+1}$ is onto.

It is known that
the map $[B\times I\ \mbox{rel}\ \partial\ ,G/O]\to
[B\times I\ \mbox{rel}\ \partial\ ,G/Top]$
has a finite cokernel (see e.g.~\cite[page 213]{Dav}), and hence so does
the map $[B\times I\ \mbox{rel}\ \partial\ ,G/O] \to
L_{n+k+1}(\mathbb Z\pi)$.

By exactness of the smooth surgery exact sequence, the
$L_{n+k+1}(\mathbb Z\pi)$-action on
$S_O(B)$ has finite orbits.

Since $h$ preserves the total Pontrjagin class, $h^m$ is tangential
for some $m>0$, so replacing $h$ by $h^m$, we can assume that
$h$ is tangential. Hence for any integer $l>0$, $h^l$ is tangential
so that the normal invariant of $[B,h^l]\in S_O(B)$ lies in
the image of the map $[B,SG]\to [B, G/Top]$ induced
by the fibration $STop\to SG\to G/Top$.
Since $SG$ is rationally contractible, $[B,SG]$ is a finite set,
so there exists an infinite sequence of positive integers $l_k$
such that the elements $[B,h^{l_k}]\in S_O(B)$
have the same normal invariant.
By exactness, $[B,h^{l_k}]$ lie in the same
$L_{n+k+1}(\mathbb Z\pi)$-orbit which is a finite set by above.
In particular, for some $p>q>0$ we have
$[B,h^{p}]=[B,h^{q}]$. Thus, for some self-diffeomorphism
$f$ of $B$, we get that $fh^{q}$ and $h^{p}$ are homotopic, or
$f$ is homotopic to $h^{p-q}$,
as wanted.
\end{proof}

\section{Vector bundles with prescribed characteristic
classes}
The following lemma is probably well-known, yet there seems to be
no reference available, so we include a complete proof.
\begin{lem}\label{lem: bundles}
Let $X$ be a finite CW-complex, $n$ be a positive integer.\newline
(i) If $k=2n+1$, then for any $n$-tuple $(p_1,\dots, p_n)$
of cohomology classes with $p_{i}\in H^{4i}(X)$ for $i=1,\dots, n$,
there is an integer $m>0$ and an orientable rank
$k$ vector bundle $\xi$ over $X$ such that
$p_{i}(\xi)=mp_{i}$ for $i=1,\dots, n$. \newline
(ii) If $k=2n$, then for any $n$-tuple $(p_1,\dots, p_{n-1},e)$
of cohomology classes with $p_{i}\in H^{4i}(X)$ for $i=1,\dots, n-1$,
and $e\in H^{2n}(X)$,
there is an integer $m>0$ and an orientable rank
$k$ vector bundle $\xi$ over $X$ such that
$e(\xi)=me$, $p_{i}(\xi)=mp_{i}$ for $i=1,\dots, n-1$.
\end{lem}
\begin{proof}
We only give a proof for $k=2n+1$; the even case is similar.
Let $\gamma^{k}$ be the universal $k$-bundle over $BSO(k)$.
We think of $p_i(\gamma^k)\in H^{4i}(BSO(k))\cong [BSO(k),K(\mathbb Z,4i)]$
as a map $BSO(k)\to K(\mathbb Z,4i)$.
It is well-known that the map
\[c=(p_1(\gamma^k),\dots,p_n(\gamma^k))\co BSO(k)\to K=
K(\mathbb Z,4)\times\dots\times K(\mathbb Z,4n)\]
is a rational homotopy equivalence and thus the homotopy groups of its
homotopy fiber $F$ are torsion. Moreover since $BSO(k)$ and $K(Z,m)$ admit a CW structure with finitely many cells in every dimension, all homotopy groups of $F$ are finitely generated.

Similarly, consider the map $f=(p_{1},\ldots,p_n)\co X \to K$
and try to lift it to $BSO(k)$. In other words, try to find
$\tilde{f}\co X\to BSO(k)$ which would make the following diagram commute
up to homotopy:
\[
\xymatrix{F\ar[r]&BSO(k)\ar[r]^{c}&K\\
&&X\ar[u]^{f}\ar@{-->}[ul]^{\tilde{f}}}
\]

We first study the auxiliary problem of trying to lift the identity map $g=\mathrm{id}\co K\to K$ to a map $\tilde g\co K\to BSO(k)$

\[
\xymatrix{F\ar[r]&BSO(k)\ar[r]^{c}&K\\
&&K\ar[u]^{g}\ar@{-->}[ul]^{\tilde{g}}}
\]

The obstructions to lifting $g$
lie in finite  groups $H^{*+1}(K;\pi_{*}(F))$, and are
generally nonzero.

Each factor $K(\mathbb Z,i)$ of $K$ is an $H$-space, and hence $K$ is too.  Let
$\times_i^{m}\co K(\mathbb Z,i)\to K(\mathbb Z,i)$ (respectively  $\times^{m}\co K\to K$)  be the
the $m$th power map.

It is easy to see that $(\times_i^{m})^*\co H^i(K(\mathbb Z,i);\mathbb Z)\to H^i(K(\mathbb Z,i);\mathbb Z)$ is multiplication by $m$.

By \cite[Lemma 4]{HQ17} the induced map $(\times_i^{m})^*\co H^{>0}(K(\mathbb Z, i);\mathbb Z/m)\to H^{>0}(K(\mathbb Z, i);\mathbb Z/m)$ is identically zero for any $m>1, i>0$.

Let $o_{j}(g)\in H^{j+1}(K;\pi_{j}(F))$ be the first nontrivial obstruction to lifting $g$. Since $\pi_{j}(F)$ is finite, by above we can find an $m$ such that $(\times^{m})^*(o_j(g))=0$.

Naturality of obstructions under pullbacks gives $o_{j}(\times^{m})=o_{j}(g\circ \times^{m})= (\times_i^{m})^*(o_j(g))=0$.

Let $l$ be the largest dimension of a cell in $X$.
By  repeating the above process finitely many times we find some $m=m(l)$ such
that all the obstructions to lifting $\times^{m}$  on the $l$-skeleton $K_l$ of $K$ vanish. Since $f\co X\to K$ can be homotoped to have the image in $K_l$ the map $\times^{m}\circ f$ can be lifted as well.

Thus, there is a map $\tilde f\co X\to BSO(k)$ such that
$c\circ\tilde f$ is homotopic to $\times^{m}\circ f$.
\[
\xymatrix{
&&BSO(n)\ar[d]^{c}\\
X\ar[r]_{f}\ar[urr]^{\tilde f}&K_l\ar[ur]_{\tilde \times^{m}}\ar[r]_{\times^{m}}&K
}
\]
Now the bundle $\xi={\tilde f}^{\#}(\gamma^{k})$ has the desired
properties.
\end{proof}

\small
\bibliographystyle{amsalpha}

\

Igor Belegradek\\School of Mathematics\\ Georgia Institute of
Technology\\Atlanta, GA, USA 30332-0160\

{\normalsize
{\it email:} \texttt{ib@math.gatech.edu}}

\

Vitali Kapovitch\\
Department of Mathematics, University of Toronto\\
Toronto, ON, Canada M5S 2E4

{\normalsize
{\it email:} \texttt{vtk@math.toronto.edu}}

\end{document}